\def\no{\if01}
\def\iftwelvept{\no}

\def\ifusepdf{\no}
\def\ifpsfont{\no}

\iftwelvept
\documentclass[leqno,12pt]{amsart}
\else
\documentclass[leqno,10pt]{amsart}
\fi
\usepackage{amssymb}
\usepackage{amscd}
\usepackage{latexsym}
\usepackage{verbatim}
\usepackage[all]{xy}

\usepackage{thmtools}

\ifusepdf
\usepackage{hyperref}
\else\fi
\ifpsfont
\usepackage[T1]{fontenc}
\usepackage{times}
\else\fi

\iftwelvept
\else\fi

\theoremstyle{plain}
\newtheorem{Theorem}{Theorem}[section]

\newtheorem{Proposition}[Theorem]{Proposition}
\newtheorem{Lemma}[Theorem]{Lemma}
\newtheorem{Corollary}[Theorem]{Corollary}

\newtheorem{Claim}{Claim}[Theorem]

\theoremstyle{definition}

\newtheorem{Definition}[Theorem]{Definition}
\newtheorem{Remark}[Theorem]{Remark}
\newtheorem{Construction}[Theorem]{Construction}
\newtheorem{Example}[Theorem]{Example}

\newtheorem{Notation}[Theorem]{Notation}

\newcommand{\ZZ}{\mathbf{Z}}

\newcommand{\DD}{\mathbb{D}}

\newcommand{\NNNN}{\operatorname{N}}

\newcommand{\AAA}{\mathcal{A}}
\newcommand{\BBB}{\mathcal{B}}

\newcommand{\DDD}{\mathcal{D}}
\newcommand{\uni}{\mathbf{1}}
\newcommand{\CCC}{\mathcal{C}}

\newcommand{\PR}{\operatorname{Pr}^{\textup{L}}}

\newcommand{\PRT}{\textup{Pr}^{\textup{t}}}
\newcommand{\PRTT}{\textup{Pr}^{\textup{t}+}}
\newcommand{\PRTTT}{\textup{Pr}^{\textup{t}\pm}}

\newcommand{\OO}{{\mathcal{O}}}

\newcommand{\MMM}{\mathcal{M}}

\newcommand{\PPP}{\mathcal{P}}

\newcommand{\Hom}{\operatorname{Hom}}

\newcommand{\Spec}{\operatorname{Spec}}

\newcommand{\Perf}{\operatorname{Perf}}

\newcommand{\SP}{\operatorname{Sp}}

\newcommand{\Mod}{\operatorname{Mod}}

\newcommand{\SSS}{\mathcal{S}}

\newcommand{\colim}{\operatorname{colim}}

\newcommand{\Cat}{\textup{Cat}_{\infty}}

\newcommand{\Map}{\operatorname{Map}}

\newcommand{\Fun}{\operatorname{Fun}}
\newcommand{\Alg}{\operatorname{Alg}}

\newcommand{\End}{\operatorname{End}}

\newcommand{\FIN}{\operatorname{Fin}_\ast}

\newcommand{\wCat}{\widehat{\textup{Cat}}_{\infty}}

\newcommand{\CAlg}{\operatorname{CAlg}}

\newcommand{\Art}{\operatorname{Art}}

\newcommand{\QC}{\operatorname{QC}}

\newcommand{\DDDD}{\mathsf{D}_{\textup{t}+}}

\newcommand{\LM}{\operatorname{\mathcal{LM}}}

\newcommand{\hhh}{\operatorname{h}}
\newcommand{\Ind}{\operatorname{Ind}}

\newcommand{\Coh}{\operatorname{Coh}}

\newcommand{\assoc}{\operatorname{As}}

\newcommand{\eone}{\mathbf{E}_1}
\newcommand{\etwo}{\mathbf{E}_2}
\newcommand{\eenu}{\mathbf{E}_n}
\newcommand{\einf}{\mathbf{E}_\infty}

\newcommand{\LMod}{\operatorname{LMod}}
\newcommand{\RMod}{\operatorname{RMod}}

\newcommand{\HL}{\widehat{L}}

\newcommand{\Proof}{{\sl Proof.}\quad}
\newcommand{\QED}{{\unskip\nobreak\hfil\penalty50\quad\null\nobreak\hfil
{$\Box$}\parfillskip0pt\finalhyphendemerits0\par\medskip}}

\begin{document}

\title{Categorified Koszul duality of algebras}

\author{Isamu Iwanari}

\begin{abstract}
Koszul duality is a duality between algebras that provides deep connections between seemingly different algebraic objects and has important applications in various areas of mathematics. The classical form of Koszul duality concerns associative algebras.
In this paper we develop a categorified generalization of Koszul duality that treats duality phenomena among monoidal stable $\infty$-categories. In our framework, monoidal stable $\infty$-categories play the role of associative algebras in the classical theory.
We establish Koszul duality results for module $\infty$-categories associated with Artin algebras and related algebras over the little 2-discs operad.
The resulting duality exhibits several structural features, including connections to right and left complete 
$t$-structures and to categories of Ind-coherent modules.

\end{abstract}






\address{Mathematical Institute, Tohoku University, Sendai, Miyagi, 980-8578, Japan}

\email{isamu.iwanari.a2@tohoku.ac.jp}

\maketitle

\section{Introduction}

Koszul duality is a duality between algebras that provides deep connections between seemingly different algebraic objects and has nontrivial applications in various areas of mathematics.  
Let \( k \) be a base field, and let \( A \) be an associative \( k \)-algebra equipped with an augmentation \( A \to k \).  
The Koszul dual \( \DD_1(A) \) is defined as the (derived) endomorphism algebra \( \End_A(k) \) of the \( A \)-module \( k \), which naturally inherits an augmentation
\[
\DD_1(A) = \mathrm{End}_A(k) \to \mathrm{End}_k(k) = k.
\]
In favorable cases, the biduality map \( A \to \DD_1(\DD_1(A)) \) is an equivalence, and \( A \) and \( \DD_1(A) \) exhibit rich interrelations.  
For instance, there exists a Fourier--Morita-type equivalence between \( A \)-modules and \( \DD_1(A) \)-modules.

\medskip

\noindent
{\it Koszul duality in the context of presentable stable \( \infty \)-categories.}
In this paper, we study a \emph{categorified} generalization of Koszul duality.  
We consider a monoidal \( \infty \)-category \( \mathcal{A} \) endowed with a monoidal functor
\[
\mathcal{A} \to \Mod_k,
\]
where \( \Mod_k \) denotes the stable \( \infty \)-category of \( k \)-module spectra (i.e., the derived \( \infty \)-category of \( k \)-vector spaces).  
We view this functor as a categorified analog of the augmentation \( A \to k \).  
In this setting, the Koszul dual \( \mathsf{D}(\mathcal{A}) \) is defined to be the endomorphism monoidal \( \infty \)-category
\[
\mathsf{D}(\mathcal{A}) = \End_{\mathcal{A}}(\Mod_k)
\]
which classifies $\AAA$-module functors $\Mod_k\to\Mod_k$. The monoidal structure is induced by composition. 
We focus on the case where \( \mathcal{A} \) is the stable \( \infty \)-category of modules over an (Artin or complete) algebra---arguably the first nontrivial instance.  
Let \( A \to k \) be an augmented \( \etwo \)-algebra over \( k \), and let
$\LMod_A \to \Mod_k$
be the induced monoidal functor given by base change.  
We then define
\[
\mathsf{D}(\LMod_A) = \End_{\LMod_A}(\Mod_k).
\]
Applying \( \mathsf{D} \) twice yields a biduality functor
\[
\delta_{\LMod_A} : 
\LMod_A \to 
\mathrm{End}_{\mathrm{End}_{\LMod_A}(\Mod_k)}(\Mod_k)
= \mathsf{D}(\mathsf{D}(\LMod_A)).
\]

Our first main result is as follows (see Theorem~\ref{main0} for details).

\begin{Theorem}
\label{intro0}
Suppose that \( A \) is Artin. Then the biduality map
\[
\delta_{\LMod_A} : \LMod_A \to \mathsf{D}(\mathsf{D}(\LMod_A))
\]
is an equivalence after passing to the right completion of \( \mathsf{D}(\mathsf{D}(\LMod_A)) \).
In particular, if we write $\mathsf{D}_{\textup{t}+}$ for $\mathsf{D}$ followed by the right completion (cf. Notation~\ref{compD}),
it induces a canonical equivalence 
\[
\LMod_A \simeq \mathsf{D}_{\textup{t}+}(\mathsf{D}_{\textup{t}+}(\LMod_A)).
\]
\end{Theorem}

This can be regarded as a categorical generalization of 
the above $A \simeq \DD_1(\DD_1(A))$.

The Koszul dual \( \mathsf{D}(\LMod_A) \) is closely related to the \( \etwo \)-Koszul dual \( \DD_2(A) \) of \( A \).  
 \( \mathsf{D}(\LMod_A) \) can be interpreted as the left completion of the monoidal \( \infty \)-category \( \LMod_{\DD_2(A)} \) of left modules.  
Our second main result is as follows (see Theorem~\ref{main1} for details).

\begin{Theorem}
\label{intro1}
Suppose that \( A \) is Artin. There is a canonical equivalence of monoidal \( \infty \)-categories
\[
\mathsf{D}(\LMod_{\DD_2(A)}) \simeq \mathrm{Ind}(\mathrm{LCoh}(A)).
\]
\end{Theorem}

The \( \infty \)-category \( \Ind(\textup{LCoh}(A)) \) consists of Ind-coherent sheaves (or modules) on \( A \) (see, e.g., \cite{Ind}).  
When \( A \) is commutative, the symmetric monoidal \( \infty \)-category \( \Ind(\textup{Coh}(A)) \) is described in \cite{Ind}.  
Thus, this result identifies the monoidal \( \infty \)-category of Ind-coherent sheaves as the \emph{categorified Koszul dual} of the \( \infty \)-category of ordinary left \( \DD_2(A) \)-modules. 
This equivalence provides an interpretation of the monoidal structure $\mathrm{Ind}(\mathrm{LCoh}(A))$
as composition of functors $\Mod_k\to \Mod_k$ (note that $\mathsf{D}(\LMod_{\DD_2(A)})=\End_{\LMod_{\DD_2(A)}}(\Mod_k)$).
This result also reveals that the categorified Koszul duality has a more subtle and interesting feature than the duality between $A$ and $\DD_2(A)$ since $\mathsf{D}(\LMod_{\DD_2(A)})$ is not $\LMod_A$ but $\mathrm{Ind}(\mathrm{LCoh}(A))$.

Next we start with $\LMod_{\DD_2(A)}$. 

\begin{Theorem}
\label{intro3}
There exists an equivalence
\[
\widehat{L}(\LMod_{\DD_2(A)})\simeq \DDDD(\DDDD(\LMod_{\DD_2(A)}))
\]
of monoidal $\infty$-categories,
where $\widehat{L}(\LMod_{\DD_2(A)})$ is the left completion of $\LMod_{\DD_2(A)}$ with respect to the standard $t$-structure.
\end{Theorem}

\begin{Remark}
Indeed, both sides are equivalent to $\mathsf{D}(\LMod_A)$ (see Proposition~\ref{dualcompletion}).
The standard $t$-structure on $\LMod_{\DD_2(A)}$ is determined by the ``connective part" $\LMod_{\DD_2(A)}^{\le0}$,
 which is the smallest subcategory containing $\DD_2(A)$ and being closed under extensions and small colimits 
(see Remark~\ref{unifiedtstr}).
This $t$-structure is generally not left complete.
\end{Remark}

\begin{Remark}
The monoidal $\infty$-category $\LMod_{\DD_2(A)}$ can also be recovered as the algebraization of $\DDDD(\DDDD(\LMod_{\DD_2(A)}))$ (see Lemma~\ref{Indend} and Definition~\ref{algdef}). Put another way, the underlying category $\LMod_{\DD_2(A)}$ may be viewed as the regularization of $\DDDD(\DDDD(\LMod_{\DD_2(A)}))$
(see \cite{PR} for the regularization).
\end{Remark}

Our results can be understood within the general paradigm of categorifying modules and sheaves in algebraic geometry.
Here, categorification refers to the process of replacing modules and sheaves by categories and sheaves of categories (cf. \cite{Ga}, \cite{TV}).
From this perspective, the categorified Koszul duality may be regarded as a categorification of the classical Koszul duality for modules.
This line of thought is partly motivated by the close relationship with categorical invariants of Landau–Ginzburg models, including matrix factorizations \cite{AL}, as well as by deformation theory of stable $\infty$-categories \cite{HD}.
For instance, let us view results from deformations of stable $\infty$-categories.
Theorem~\ref{intro1} yields a functor
\[
\RMod_{\LMod_{\DD_2(A)}}(\PRTT)\longrightarrow \LMod_{\mathsf{D}(\LMod_{\DD_2(R)})}(\PRT)\simeq \LMod_{\Ind(\textup{LCoh}(A))}(\PRT)
\]
which carries $\mathcal{E}$ to $\Fun_{\LMod_{\DD_2(A)}}(\Mod_k, \mathcal{E})$.
Here $\PRT$ (resp. $\PRTT$) denote the symmetric monoidal $\infty$-category of presentable $\infty$-categories with accessible $t$-structures
(resp. with right complete accessible $t$-structures),
see Section~\ref{tsec}.
Using Proposition~\ref{pullbackdual}, we have $\widehat{R}_A:\LMod_{\Ind(\textup{LCoh}(A))}(\PRT)\to \LMod_{\LMod_A}(\PRTT)$ induced
by the right completion functor $\widehat{R}:\PRT\to \PRTT$.
Composing these, we obtain the functor
\[
\kappa:\RMod_{\LMod_{\DD_2(A)}}(\PRTT)\longrightarrow \LMod_{\LMod_A}(\PRTT)
\]
given informally by $\widehat{R}(\Fun_{\LMod_{\DD_2(A)}}(\Mod_k,-))$.
This is a {\it right adjoint} to
\[
\flat:\LMod_{\LMod_A}(\PRTT)   \longrightarrow  \RMod_{\LMod_{\DD_2(A)}}(\PRTT)
\]
which is given by the base change $\DDD\mapsto \Mod_k\otimes_{\LMod_A}\DDD$.
Furthermore, Theorem~\ref{intro1} and Proposition~\ref{pullbackdual} imply that $\kappa$ is fully faithful as a $\PRTT$-enriched functor when restricted to the full subcategory
of objects generated by $\Mod_k$ under small limits (e.g., stable $\infty$-categories of representations of groups).
Since $\flat$ realizes the base change of $R$-linear stable presentable $\infty$-categories with right complete $t$-structures from $R$ to $k$, the right adjoint $\kappa$ plays a crucial role that associates to a class of Hochschild cohomology a deformation of stable $\infty$-categories with right complete $t$-structures.
This means that one can apply the results to deformation theory (deformations include possibly non-left complete $t$-structures), which is different from \cite{HD}.
We will discuss this application in future work.

\vspace{3mm}

This paper is organized as follows.
In Section~\ref{Pre}, we fix some notation and terminology
and review some results which we will use in this paper.
In Section~\ref{ICD}, we formulate Theorem~\ref{intro1}.
In Section~\ref{DS} we formulate Theorem~\ref{intro0} and Theorem~\ref{intro3}.
This section also includes the review of $t$-structures.
In Section~\ref{proofsect}, we complete the proofs of Theorem~\ref{intro1}, Theorem~\ref{intro0} and Theorem~\ref{intro3}.
In Section~\ref{FS} we discuss a generalization to formal stacks.  
In Appendix~\ref{app}, we provide a covariantly functorial assignment
of categories of Ind-coherent sheaves on Artin $\eone$-algebras.

\section{Preliminaries}

\label{Pre}

Let $k$ denote a fixed field, which serves as the base field.

\subsection{}

\label{notation}

{\it Notation and Terminology.}
We use 
the language of $(\infty,1)$-categories.
We use the theory of {\it quasi-categories} as a model of $(\infty,1)$-categories.
Our main references are \cite{HTT}
 and \cite{HA}.
Following \cite{HTT}, we shall refer to quasi-categories
as {\it $\infty$-categories}.
To an ordinary category, we can assign an $\infty$-category by taking
its nerve, and therefore
when we treat ordinary categories, we often omit the nerve $\NNNN(-)$
and directly regard them as $\infty$-categories.
Here is a list of some of the notation and conventions we frequently use:

\begin{itemize}

\item $\simeq$: categorical equivalence or homotopy equivalence.  The symbol $\simeq$ is used to indicate isomorphisms and equivalences in categories and
$\infty$-categories.

\item $\CCC^{op}$: the opposite $\infty$-category of an $\infty$-category.

\item $\Map_{\mathcal{C}}(C,C')$: the mapping space from an object $C\in\mathcal{C}$ to $C'\in \mathcal{C}$ where $\mathcal{C}$ is an $\infty$-category.

\item $\Fun(A,B)$: the function complex for simplicial sets $A$ and $B$. If $A$ and $B$ are $\infty$-categories, we regard $\Fun(A,B)$ as the functor category.

\item $\Cat$: the $\infty$-category of small $\infty$-categories. We write $\wCat$ for the $\infty$-category of $\infty$-categories in an enlarged universe. 

\item $\Ind(\CCC)$: the Ind-category of $\CCC$.

\end{itemize}

We will use the theory of higher algebras.
We refer the reader to \cite{HA}.
Here is a list of (some) of the notation about $\infty$-operads and algebras over them that we will use:

\begin{itemize}
\label{notation2}

\item $\Gamma$: $\FIN$  denotes the category of pointed finite
sets $\langle 0 \rangle, \langle 1 \rangle,\ldots \langle n \rangle,...$
where $\langle n \rangle=\{*,1,\ldots n \}$ with the base point $*$.
We write $\Gamma$ for $\NNNN(\FIN)$. $\langle n\rangle^\circ=\langle n\rangle\backslash*$.

\item $\Alg_n(\MMM)$:
For a symmetric monoidal $\infty$-category $\MMM^\otimes$,
we write $\Alg_{n}(\MMM)$ for
the $\infty$-category of algebra objects over the $\infty$-operad of little $n$-cubes $\mathbf{E}^\otimes_n$ in $\MMM^\otimes$  (cf. \cite[5.1]{HA}).
We refer to an object of $\Alg_{n}(\MMM)$ as an $\eenu$-algebra in $\MMM$.
If we denote by $\assoc^\otimes$ the associative operad (\cite[4.1.1]{HA}),
there is the standard equivalence $\assoc^\otimes\simeq \eone^\otimes$
of $\infty$-operads. We usually identify $\Alg_{1}(\MMM)$
with the $\infty$-category $\Alg_{\assoc}(\MMM)$, that is, the 
$\infty$-category of associative algebras in $\MMM$.

\item $\CAlg(\mathcal{M})$: the $\infty$-category of commutative
algebra objects in a symmetric
monoidal $\infty$-category $\mathcal{M}^\otimes$.
Instead of $\Alg_{\einf}(\MMM)$ or $\Alg_{\textup{comm}}(\MMM^\otimes)$
we write $\CAlg(\mathcal{M})$ for the $\infty$-category of algebra objects over the $\einf$-operad or the commutative operad in $\MMM^\otimes$.

\item $\Mod_R(\mathcal{M})$: 
the $\infty$-category of
$R$-module objects,
where $\mathcal{M}^\otimes$
is a symmetric monoidal $\infty$-category and $R\in \CAlg(\MMM)$.

\item $\Mod_R$: When $\MMM^\otimes$ is the symmetric monoidal $\infty$-category $\SP^\otimes$
of spectra, 
we write $\Mod_R^\otimes$
for
$\Mod_R^\otimes(\SP)$.
We denote by $\Mod_R$ the underlying category, that is, $\Mod_R(\SP)$.
When $R$ is the Eilenberg-Maclane ring spectrum $HA$, we write $\Mod_A$ for $\Mod_{HA}$.
For example, $\Mod_k$ is equivalent to the derived $\infty$-category of $k$-modules.

\item $\CAlg_k$: $\CAlg_k=\CAlg(\Mod_k)$.

\item 
$\LMod(\MMM)$:
the $\infty$-category $\Alg_{\LM}(\MMM)$ of $\LM^\otimes$-algebras
over $\LM^\otimes$: the $\infty$-operad defined in \cite[4.2.1.7]{HA}.
We may consider an object of  $\LMod(\MMM)$ to be a pair 
of $(A,M)$ such that $A$ is an associative algebra object in $\MMM$
and a left $A$-module $M$.
There is a Cartesian fibration $\LMod(\MMM)\to \Alg_{\assoc}(\MMM)\simeq \Alg_1(\MMM)$
which sends $(A,M)$ to $A$.
Another canonical morphism is the forgetful functor $\LMod(\MMM)\to \MMM$
which is informally given by $(A,M)\mapsto M$.

\item $\LMod_A(\MMM)$, $\LMod_A$: the $\infty$-category of left $A$-module objects in $\MMM$.
For $A\in \Alg_1(\MMM)$, we define $\LMod_A(\MMM)$
to be the fiber of $\LMod(\MMM)\to \Alg_1(\MMM)$ over $A$ in $\Cat$.
If $\MMM$ is the $\infty$-category of spectra (or $\Mod_k$),
we write $\LMod_A$ for $\LMod_A(\MMM)$.

\item $\RMod(\MMM)$: the right module version of $\LMod(\MMM)$

\item $\RMod_B(\MMM)$, $\RMod_A$: the $\infty$-category of right $B$-module objects in $\MMM$.
$\RMod_A$ is the right module version of $\LMod_A$.

\item $\PR$: $\infty$-category of presentable $\infty$-categories such that the mapping space consists of functors which preserve small colimits. 

\item $\PR_{\mathbb{S}}$, $\PR_k$: $\PR_{\mathbb{S}}$ is the $\infty$-category of presentable stable $\infty$-categories. A morphism $\CCC\to \DDD$ in $\PR_{\mathbb{S}}$ is a functor that preserves small colimits. It has a closed symmetric monoidal structure whose unit is the stable $\infty$-category $\SP$ of spectra  (see \cite[Definition 5.5.3.1]{HTT}, \cite[4.8]{HA}).
We write $\PR_k$ for $\Mod_{\Mod_k}(\PR_{\mathbb{S}})$ and write $\otimes_k$ for the tensor products.

\end{itemize}

We define the functor categories.
Let $\AAA^\otimes$ be a $k$-linear monoidal stable $\infty$-category i.e., $\AAA^\otimes \in \Alg_1(\PR_k)$.
We write $\PR_{\AAA}:=\LMod_{\AAA^\otimes}(\PR_k)$.
We note that $\PR_{\AAA}$ is right-tensored over $\PR_{k}$ (we denote the tensor by $\otimes_k$).
Let $\BBB,\CCC \in \PR_{\AAA}$. 
We denote by $\Fun_{\AAA}(\BBB,\CCC)\in \PR_k$ the Hom object which has the following universal property:
There exists a morphism $\mathcal{B}\otimes_k \Fun_{\AAA}(\BBB,\CCC)\to \CCC$ in $\PR_{\AAA}$ such that it determines
an equivalence of spaces
\[
\Map_{\PR_k}(\mathcal{E}, \Fun_{\AAA}(\BBB,\CCC))\simeq \Map_{\PR_{\AAA}}(\BBB\otimes_k \mathcal{E},\CCC)
\]
for any $\mathcal{E}\in \PR_k$.
The existence follows from a standard argument, see e.g., \cite[Lemma 5.1]{IH}.
When $\BBB=\CCC$, we write $\End_{\AAA}(\BBB)=\Fun_{\AAA}(\BBB,\CCC)$.
In this case it is promoted to an object of $\Alg_1(\PR_k)$ (see \cite[4.7.1.40, 4.7.1.41]{HA}).
We refer to $\End_{\AAA}(\BBB)$ as the endomorphism algebra.
For $R\in \CAlg(\SP)$, we often write $\Fun_R(-,-)$ for $\Fun_{\Mod_R}(-,-)$.

\subsection{}
\label{koszulsec}
Let $\MMM^\otimes$ be a symmetric monoidal $\infty$-category which has small colimits and limits such that the tensor product preserves sifted colimits. We write $\uni$ for a unit object.
We review the Koszul duality of algebras in $\MMM^\otimes$ in the formalism of pairings (see \cite{HA}, \cite[X]{DAG} for more details).
Let $\Alg^+_{1}(\MMM)=\Alg_1(\MMM)_{/\uni}$ be the $\infty$-category of augmented $\eone$-algebra objects (associative algebra objects) in $\MMM^\otimes$.
For a pair $(\epsilon_A:A\to \uni,\epsilon_B:B \to \uni)\in \Alg^+_{1}(\MMM)\times \Alg^+_{1}(\MMM)$ we put
\[
\textup{Pair}(\epsilon_A,\epsilon_B):= \Map_{\Alg_{1}(\MMM)}(A\otimes B,\uni)\times_{\Map(A,\uni)\times\Map(B,\uni)}\{(\epsilon_A,\epsilon_B)\}.
\]
Fix $\epsilon_B:B \to \uni$ and consider the functor
$u_{\epsilon_B}:\Alg_1^+(\MMM)^{op}\to \SSS$ given informally by $[\epsilon_A:A\to \uni]\mapsto \textup{Pair}(\epsilon_A,\epsilon_B)$. 
The Koszul dual $\DD_{\MMM^\otimes}(B)\in \Alg_1^+(\MMM)$ of $B$ (we abuse notation since $\DD_{\MMM^\otimes}(B)$ depends on the augmentation $\epsilon_B:B \to \uni$)
is defined by the following universal property that $\DD_{\MMM^\otimes}(B)$ represents the functor $u_{\epsilon_B}$.
Under a mild condition on $\MMM^\otimes$, 
there exist $\DD_{\MMM^\otimes}(B)$ and functorial equivalences
\[
\Map_{\Alg_1^+(\MMM)}(A,\DD_{\MMM^\otimes}(B))\simeq \textup{Pair}(\epsilon_A,\epsilon_B)
\]
for $\epsilon_A:A\to \uni$. We shall refer to the pairing $B\otimes\DD_{\MMM^\otimes}(B)\to \uni$ 
corresponding to the identity of $\DD_{\MMM^\otimes}(B)$ as the universal pairing for $B\to k$.
In that case, there is a functor 
\[
\DD_{\MMM^\otimes}:\Alg_1(\MMM) \longrightarrow \Alg_1(\MMM)^{op}
\]
which carries $A\to \uni$ to $\DD_{\MMM^\otimes}(A)\to \uni$.

The Koszul duals can be described in terms of bar constructions and endomorphism algebras.
Let $\textup{Bar}:\Alg_{1}^+(\MMM)\to \MMM$ be the functor
given by $[B\to \uni] \mapsto \uni\otimes_B\uni$, which we refer to as the bar construction of augmented algebras.
Then $\uni\otimes_B\uni$ admits a structure of a coaugmented coalgebra in a suitable way, and $\textup{Bar}$
is promoted to $\textup{Bar}:\Alg_{1}^+(\MMM)\to \Alg^+_{1}((\MMM)^{op})^{op}$
(notice that the target is the $\infty$-category of coaugmented coalgebras).
According to \cite[X, Example 4.4.20]{DAG}, for $B\in \Alg^+_{1}(\MMM)$, the Koszul dual $\DD_{\MMM^\otimes}(B)$
is equivalent to
the dual $\textup{Bar}(B)^\vee$ in $\MMM$, that is, the endomorphism algebra $\End_{\LMod_B(\MMM)}(\uni)$.

We recall the Koszul duality between $\eenu$-algebras. In this note, we assume that $\MMM^\otimes =\Mod_k$: in this case
Koszul duals exist (cf. \cite[X, 3.1.5, 4.3.8, 4.4.8]{DAG}).
Let $n\ge1$ be a natural number and  let $\Alg_{n}(\Mod_k)$ be the $\infty$-category of $\eenu$-algebras in $\Mod_k$.
We write $\Alg^+_{n}(\Mod_k)=\Alg_n(\Mod_k)_{/k}$ for the $\infty$-category
$\Alg_{n}(\Mod_k)_{/k}$ of augmented $\eenu$-algebras.
We consider the Koszul duals of augmented $\eenu$-algebras.
We take $\MMM^\otimes$ to be $\Alg_{n-1}^\otimes(\Mod_k)$. There exists the Koszul duality functor
\[
\DD_n:\Alg_{n}^+(\Mod_k)\simeq \Alg_1^+(\Alg_{n-1}(\Mod_k))\longrightarrow\Alg_1^+(\Alg_{n-1}(\Mod_k))^{op}\simeq  \Alg_{n}^+(\Mod_k)^{op}
\]
(see e.g. \cite[X, 4.4.8]{DAG}).
Let $\epsilon_{B}:B\to k$ be an augmented $\eenu$-algebra.
As mentioned above,
$\DD_n(B)$ is determined by the following universal property.
For an augmented $\eenu$-algebra $\epsilon_C:C\to k$, we write
\[
\textup{Pair}(\epsilon_B,\epsilon_C)= \Map_{\Alg_{n}(\Mod_k)}(B\otimes_k C,k)\times_{\Map(B,k)\times\Map(C,k)}\{(\epsilon_B,\epsilon_C)\}.
\]
Then the functor $\Alg^+_{n}(\Mod_k)^{op}\to \SSS$ given by $[\epsilon_C:C\to k]\mapsto \textup{Pair}(\epsilon_B,\epsilon_C)$
is representable by $\DD_{n}(B)$.
There is a universal (Koszul dual) pairing $B\otimes_k\DD_{n}(B)\to k$ corresponding to
the identity map $\textup{id}\in \Map_{\Alg^+_{n}(\Mod_k)}(\DD_{n}(B),\DD_{n}(B))$.
We refer to $\DD_{n}(B)$ as the $\eenu$-Koszul dual of $B$.


Next, we consider the categorified context, that is, the case of $\MMM^\otimes=\PR_k$. 
There is a self-dual adjoint functor  
\[
\xymatrix{
 \mathsf{D}:=\DD_{(\PR_k)^\otimes}: \Alg_1(\PR_k)_{/\Mod_k^\otimes} \ar@<0.5ex>[r] & (\Alg_1(\PR_k)_{/\Mod_k^\otimes})^{op}
}
\]
which we shall refer to as the categorified Koszul dual functor. 
The functor $\mathsf{D}$ sends $\mathcal{A}^\otimes\to \Mod_k^\otimes$ to the $\Mod_k^\otimes$-linear dual
of the bar construction $\Mod_k\otimes_{\mathcal{A}}\Mod_k$, that is, the endomorphism algebra
object $\End_{\mathcal{A}}(\Mod_k)\in \Alg_1(\PR_k)$
of $\Mod_k$ in $\PR_{\mathcal{A}}$, that is equipped with the augmentation $\End_{\mathcal{A}}(\Mod_k)\to \End_{\Mod_k}(\Mod_k)\simeq \Mod_k$
(see \cite[X, Example 4.4.20]{DAG}).

\subsection{Artin algebras}

\begin{Definition}
Let $A$ be an augmented $\eenu$-algebra in $\Mod_k$. Let $n$ be a natural number or $\infty$.
The augmented $\eenu$-algebra $R$ is Artin if the following properties hold:
\begin{itemize}
\item $A$ is connective,

\item $\oplus_{i\in \ZZ}H^{i}(A)$ is finite dimensional over $k$,

\item if $I$ is the radical of $H^0(A)$, the canonical map $k\to H^0(A)/I$ is an isomorphism.

\end{itemize}
By convention, an $\mathbf{E}_{\infty}$-algebra is a commutative algebra object in $\Mod_k$, and $\Alg_\infty(\Mod_k)=\CAlg(\Mod_k)$.
We let $\textup{Art}_n$ denote the full subcategory on $\Alg_n^+(\Mod_k)$ spanned by Artin $\eenu$-algebras.
When $n=\infty$, we write $\textup{Art}=\textup{Art}_n$.
\end{Definition}

In \cite[X]{DAG}, the terminology small is adopted instead of Artin.

\begin{Example}
For $s\ge 0$, the trivial square zero extension $k\oplus k[s]$ of $k$ by $k[s]$ is Artin.
\end{Example}

\begin{Remark}
\label{augnonaug}
Let $\Art_n^\circ$ be the full subcategory of $\Alg_n(\Mod_k)$ spanned by the essential image of the functor $U:\Art_n\to \Alg_n(\Mod_k)$
determined by forgetting augmentations.
Then it gives an equivalence $\Art_n\simeq \Art_n^\circ$.
We see this by observing that for any $A\in \Art_n$ the mapping space $\Map_{\Alg_n(\Mod_k)}(U(A),k)\simeq \Map_{\Alg_n(\Mod_k)}(H^0(U(A)),k)$ is contractible.
\end{Remark}

\begin{Remark}
\label{smallrem}
For $A\in \textup{Art}_n$
there exists a sequence
\[
A=A(0)\to A(1)\to \cdots \to A(m)=k
\]
in $\textup{Art}_n$
such that,
for each $0\le i\le m-1$,
there exists a pullback square
\[
\xymatrix{
A(i) \ar[r] \ar[d] & A(i+1) \ar[d] \\
k \ar[r] &  k\oplus k[s]
}
\]
with $s\ge1$ (cf. \cite[X, Proposition 4.5.1]{DAG}).
\end{Remark}

\begin{Remark}
When $A$ is Artin, the biduality map $A\to \DD_n(\DD_n(A))$ is an equivalence.
\end{Remark}

\subsection{Ind-coherent sheaves}
\label{Indsec}
We will give a minimal review of Ind-coherent sheaves. The basic reference is \cite{Ind} (see also Appendix for the case of Artin $\eenu$-algebras).

Let $A$ be a connective $\eenu$-algebra in $\Mod_k$. Suppose that $H^0(A)$ is Artin. 
Let $\textup{LCoh}(A)\subset \LMod_A$ 
be the stable subcategory spanned by those objects that have bounded amplitudes with respect to the standard $t$-structure (cf. Example~\ref{connectivetstr})
such that each cohomology is finitely generated as an $H^0(A)$-module.
For ease of notation, we often write $\Coh(A)$ for $\textup{LCoh(A)}$.
We let $\Ind(\textup{LCoh}(A))$ be the $\infty$-category
of Ind-objects (see \cite[5.3]{HTT} for Ind-objects).
In this case, $\textup{LCoh}(A)$ can be characterized as  the smallest stable subcategory
of $\LMod_A$ which contains the residue field $k\simeq H^0(A)/I$ and is closed under retracts (cf. e.g. \cite[X, Remark 3.4.2]{DAG}).

Let $A$ and $B$ be connective $\eenu$-algebras.
Suppose that $H^0(A)$ and $H^0(B)$ are Artin, that is, $H^0(A),H^0(B)\in \Art_n$.
Let $f^\sharp:A\to B$ be a morphism.
From the definition of $\textup{LCoh}(A)$ and $\textup{LCoh}(B)$ and the assumption on $A$ and $B$,
the restriction functor $\LMod_B\to \LMod_A$ along $f^\sharp$ carries $\textup{LCoh}(B)$ to $\textup{LCoh}(A)$.
It is uniquely extended to a colimit-preserving ($\Mod_k$-linear) functor 
$f_*:\Ind(\textup{LCoh}(B))\rightarrow \Ind(\textup{LCoh}(A))$. We refer to it as the (proper) pushforward functor along $f$.
We define the $!$-pullback functor to be the right adjoint of $f_*$
\[
f^!: \Ind(\textup{LCoh}(A))\longrightarrow \Ind(\textup{LCoh}(B)).
\]
See also Lemma~\ref{dualind} (3).
This $!$-pullback functor preserves small colimits since the left adjoint $f_*$ preserves compactness of objects.

We recall the monoidal structure.
Let $A$ be a connective commutative $k$-algebra  (see \cite[5.6]{Ind}).
(See Appendix Definition~\ref{eenuindmodule} for the case of Artin $\etwo$-algebras: In that situation $\Ind(\textup{LCoh}(A))$ admits an associative monoidal structure.)
The stable $\infty$-category $\Ind(\Coh(A))$ has a symmetric monoidal structure such that the tensor product is defined as the $!$-pullback functor 
\[
\Delta^!:\Ind(\Coh(A))\otimes_k\Ind(\Coh(A))\simeq \Ind(\Coh(A\otimes_kA))\to \Ind(\Coh(A))
\]
of the multiplication map
$\Delta^\sharp:A\otimes_kA\to A$ (geometrically corresponding to the diagonal morphism)
where the first equivalence is the tensor compatibility of Ind-coherent sheaves 
(see \cite[Vol.1, Chapter 4, 6.3.3]{Gai2}). 
If we write $\omega_A$ for the dualizing object $p_A^!(k)$ where $p_A:\Spec A\to \Spec k$ is the structure morphism (see above for the $!$-pullback for the Artin case),
a unit object is $\omega_A$.
The assignment $A\mapsto \Ind(\Coh(A))$
is promoted to a symmetric monoidal functor 
\[
\theta_k^{\Ind,\infty}:\textup{Art}\longrightarrow \CAlg(\PR_k)
\]
which carries $A$ to $\Ind(\Coh(A))$ and 
$f^\sharp:A\to B$ to a morphism $f^!:\Ind(\Coh(A))\to\Ind(\Coh(B))$ in $\CAlg(\PR_k)$ (see \cite[Vol.1, Chapter 5, 4.3]{Gai2} or Appendix of this paper).
In the context of Artin $\etwo$-algebras, there is a monoidal functor 
\[
\theta_k^{\Ind,2}:\textup{Art}_2\longrightarrow \Alg_1(\PR_k)
\]
which carries $A$ to $\Ind(\textup{LCoh(A)})$
(see Definition~\ref{eenuindmodule}).

Let $A$ be an Artin $\etwo$-algebra. 
There is a monoidal fully faithful functor
\[
\Upsilon_A:\LMod_A \to \Ind(\textup{LCoh}(A))
\]
see \cite{Gai2} and Appendix.

\begin{Lemma}[cf. \cite{DAG}]
\label{koszulfunc}
The following holds.
\begin{enumerate}

\item For $A\in \textup{Art}_1$, there exists an equivalence
\[
\LMod_{\DD_1(A)}\simeq \Ind(\textup{RCoh}(A))
\]
in $\PR_k$, which is the Ind-extension of $\textup{LPerf}_{\DD_1(A)}\simeq \textup{RCoh}(A)$
given informally by $M\mapsto k\otimes_{\DD_1(A)}M$.
Similarly, there exists an equivalence
\[
\RMod_{\DD_1(A)}\simeq \Ind(\textup{LCoh}(A))
\]
in $\PR_k$, which is the Ind-extension of full subcategories of compact objects
$\textup{RPerf}_{\DD_1(A)}\simeq \textup{LCoh}(A)$
given informally by $N\mapsto N\otimes_{\DD_1(A)}k$.
Here $\textup{RPerf}_{\DD_1(A)}$ (resp. $\textup{LPerf}_{\DD_1(A)}$) is the full subcategory of $\RMod_{\DD_1(A)}$
(resp. $\textup{LPerf}_{\DD_1(A)}$), that consists of compact objects.

\item Let $f^\sharp:B\to A$ be a morphism of Artin $\eone$-algebras.
There exists a commutative diagram
\[
\xymatrix{
 \LMod_{\DD_1(A)}   \ar[r]^{\simeq}   &  \Ind(\textup{RCoh}(A)) \\
\LMod_{\DD_1(B)}  \ar[r]^{\simeq} \ar[u]^{\textup{rest}}  &  \Ind(\textup{RCoh}(B))  \ar[u]^{f^!}
}
\]
such that the horizontal equivalences come from (1), and the left vertical arrow is the restriction functor along
$\DD_1(A)\to \DD_1(B)$. Similarly, the opposite version also holds. 
\end{enumerate}
\end{Lemma}

\Proof
Indeed, (1) is proved in \cite[Theorem 3.5.1]{DAG}.
However, the constructions of equivalences are relevant to us (e.g. Lemma~\ref{standardcompletion}) so
that
we will provide the proof of (1).
The latter case is similar to the former (or it follows by taking the opposite algebra and $\DD_1(A^{op})\simeq \DD_1(A)^{op}$).
We will treat the former equivalence.
Consider 
$(-)\otimes_{\DD_1(A)}k: \textup{LPerf}_{\DD_1(A)}\to \textup{RCoh}(A)$
given informally by $M\mapsto k\otimes_{\DD_1(A)}M$.
From the description of $\DD_1(A)$ as the endomorphism algebra 
we have $\Map_{\textup{LPerf}_{\DD_1(A)}}(\DD_1(A),\DD_1(A))\simeq \Map_{\textup{RCoh}(A)}(k,k)$
so that $(-)\otimes_{\DD_1(A)}k$ is the fully faithful when the functor is restricted to the full subcategory
spanned by $\DD_1(A)$. Note that $\textup{LPerf}_{\DD_1(A)}$ is the smallest stable subcategory which contains $\DD_1(A)$ 
and is closed under retracts, and  
the stable $\infty$-category $\textup{RCoh}(A)$ itself is the smallest stable subcategory
which contains $k$ 
and is closed under retracts.
Consequently, we see that $(-)\otimes_{\DD_1(A)}k$ is an equivalence so that its Ind-extension is also equivalence. 

Next, we will prove (2). 
It is enough to show the diagram consisting of left adjoints of the desired diagram
\[
\xymatrix{
 \LMod_{\DD_1(A)} \ar[d]_{\DD_1(B)\otimes_{\DD_1(A)}}   \ar[r]^{\simeq}   &  \Ind(\textup{RCoh}(A)) \ar[d]^{f_*} \\
\LMod_{\DD_1(B)}  \ar[r]^{\simeq}  &  \Ind(\textup{RCoh}(B))  
}
\]
commutes up to homotopy.
 At the level of compact objects the diagram 
\[
\xymatrix{
 \textup{LPerf}_{\DD_1(A)} \ar[d]_{u}   \ar[r]^{s}   &  \textup{RCoh}(A) \ar[d]^{t} \\
\textup{LPerf}_{\DD_1(B)}  \ar[r]^v  &  \textup{RCoh}(B)  
}
\]
commutes
where $s$ is given by $k\otimes_{\DD_1(A)}(-)$, $t$ is the restriction functor,
$u$ is given by
$\DD_1(B)\otimes_{\DD_1(A)}(-)$, and 
$v$ is given by $(-)\otimes_{\DD_1(B)}k$.
To see this,
we note that $t\circ s$ is induced by the integral kernel that is
the $\DD_1(A)$-$B$-bimodule $k$ determined by the composition $\DD_1(A)\otimes_kB\to \DD_1(A)\otimes_kA \to k$
of the universal pairing and $\textup{id}\otimes f_*$.
On the other hand, $v \circ u$
is 
induced by the integral kernel that is
the $\DD_1(A)$-$B$-bimodule $k$ determined by the composition $\DD_1(A)\otimes_kB\to \DD_1(B)\otimes_kB \to k$
of the universal pairing and $\DD_1(f_*)\otimes \textup{id}$.
Now we complete the proof by observing that the canonically defined diagram
\[
\xymatrix{
 \DD_1(A)\otimes_kB \ar[d]^{\otimes(B\to A)} \ar[r]^{(\DD_1(A)\to \DD_1(B))\otimes}   &  \DD_1(B)\otimes_kB \ar[d]^{\textup{universal}} \\
\DD_1(A)\otimes_kA  \ar[r]^{\textup{universal}}  &  k  
}
\]
commutes.
\QED

\section{The Ind-coherent sheaves as Koszul dual}

\label{ICD}

In this section
we formulate the statement of Theorem~\ref{intro1} (see Theorem~\ref{main1}).

\subsection{}

{\it Algebraization.}
We recall from \cite[4.8.5]{HA} that
there exists a symmetric monoidal left adjoint functor 
\[
\xymatrix{
 \Theta_k:\Alg_{1}(\Mod_k) \ar@<0.5ex>[r] & (\PR_k)_{\Mod_k/}
}
\]
such that
the left adjoint $\Theta_k$ is fully faithful, and $\Theta_k(A)=[A\otimes_k(-):\Mod_k\to \LMod_A]$.
For a morphism $f^{\sharp}:A\to A'$, $\Theta_k(f)$ is the base change functor $A'\otimes_A(-):\LMod_A\to \LMod_{A'}$.
We abuse notation by denoting by
\[
\xymatrix{
 \Theta_k:\Alg_{2}(\Mod_k) \ar@<0.5ex>[r] & \Alg_1(\PR_k): E_k \ar@<0.5ex>[l]
}
\]
the induced adjunction
such that
\begin{itemize}
\item the left adjoint $\Theta_k$ is fully faithful, and $\Theta_k(A)$ is the monoidal $\infty$-category $\LMod_A^\otimes$,

\item the right adjoint $E_k$ carries $\MMM^\otimes\in \Alg_{1}(\PR_k)$ to (the opposite algebra of) the endomorphism algebra $\End(\mathsf{1}_{\MMM})$
of the unit object $\mathsf{1}_{\MMM}$, defined as an $\etwo$-algebra object in $\Mod_k$.

\end{itemize}
See e.g. \cite[Section 3.1]{IM} for endomorphism algebras of objects of presentable $\infty$-categories.
\begin{Definition}
\label{algdef}
For $\MMM^\otimes\in\Alg_{1}(\PR_k)$, we shall refer to $\Theta_k\circ E_k(\MMM^\otimes)$ as the algebraization of $\MMM^\otimes$.
We refer to the canonical monoidal functor $\Theta_k\circ E_k(\MMM^\otimes)\to\MMM^\otimes$ as the algebraization functor. 
\end{Definition}

\begin{Lemma}
\label{dualcatalg}
Let $A$ be an $\etwo$-algebra in $\Mod_k$.
The $\etwo$-Koszul dual $\DD_2(A)$ can be identified with $E_k(\End_{\LMod_{A}}(\Mod_k))\in \Alg_2(\Mod_k)$.
If we write $\kappa:\LMod_A\otimes_k\End_{\LMod_A}(\Mod_k)\to \Mod_k$ for the tautological/universal pairing, then
the universal $\etwo$-Koszul pairing $A\otimes_k\DD_2(A) \to k$
is equivalent to the composition of
\[
E_k(\kappa):E_k(\LMod_A\otimes_k\End_{\LMod_A}(\Mod_k))\to E_k(\Mod_k)=k
\]
and the canonical 
morphism $E_k(\LMod_A)\otimes_k E_k(\End_{\LMod_A}(\Mod_k))\to E_k(\LMod_A\otimes_k\End_{\LMod_A}(\Mod_k))$.
\end{Lemma}

\Proof This is proved in \cite[Lemma 3.6]{HD}.
\QED

Let $A\in \textup{Art}_2$ be an Artin $\etwo$-algebra.
Let $i^\sharp:A\to H^0(A/I)=k$ be the canonical map.
We have the adjoint pair 
\[
\xymatrix{
i_*: \Ind(\Coh(k))=\Mod_k \ar@<0.5ex>[r] & \Ind(\textup{LCoh}(A)) :i^! \ar@<0.5ex>[l]
}
\]
(cf. Section~\ref{Indsec}).
Since $A$ is Artin, $\Coh(A)=\textup{LCoh}(A)$
is the smallest stable subcategory of $\LMod_A$ which contains $k$ and is closed under retracts.
We note that $i^!$ is a monoidal colimit-preserving functor.
The underlying functor is given by the Hom complex $\Hom_{\Ind(\Coh(A))}(k,-)$ since it is a right adjoint of the evident pushforward.
The functor $i^!$ exhibits $\Mod_k$ as a left $\Ind(\Coh(A))$-module object.
We consider 
the endomorphism monoidal $\infty$-category $\End_{\Ind(\Coh(A))}(\Mod_k)$,
which is defined as an object of $\Alg_1(\PR_k)$.

Now we claim

\begin{Lemma}
\label{Indend}
Let $E_k(\mathsf{D}(\LMod_A))\simeq \DD_2(A)$ be the equivalence induced by the identification
of the $\etwo$-Koszul dual in Lemma~\ref{dualcatalg}.
Let $E_k(\End_{\Ind(\textup{LCoh}(A))}(\Mod_k))=E_k(\mathsf{D}(\Ind(\textup{LCoh})))\to E_k(\mathsf{D}(\LMod_A))$ be the morphism
obtained from $\Upsilon_A:\LMod_A\to \Ind(\textup{LCoh}(A))$ by applying $\mathsf{D}$ and $E_k$.
Then $E_k(\mathsf{D}(\Ind(\textup{LCoh})))\to E_k(\mathsf{D}(\LMod_A))\simeq \DD_2(A)$ is an equivalence.
Namely,
the endomorphism algebra of the identity functor $\Mod_k\to \Mod_k$ (i.e., the unit object) in $\End_{\Ind(\textup{LCoh}(A))}(\Mod_k)$
is equivalent to $\DD_2(A)$.

\end{Lemma}

\Proof
Indeed, according to \cite[Theorem 1.3.1]{BPN},
\[
\Phi^!:\Ind(\Coh(k\otimes_Ak))^{<\infty}\simeq \End_{\Ind(\textup{LCoh}(A))}(\Mod_k)^l
\]
where the subscript $<\infty$ means the left bounded with respect to the $t$-structure, and the superscript $l$
means the full subcategory spanned by left $t$-exact functors up to finite shifts (this proof apparently treats the case of commutative $A$, but a version of
\cite[Theorem 1.3.1]{BPN} holds for an Artin $\etwo$-algebra $A$ so that it works also for Artin $\etwo$-algebras, see Remark~\ref{etwocaserem} below).
The endomorphism algebra of $k\in \Coh(k\otimes_Ak)\subset \Ind(\Coh(k\otimes_Ak))_{<\infty}$ is $\DD_2(A)\simeq \DD_1(k\otimes_Ak)$.
This argument only says that the endomorphism algebra of the identity functor $\Mod_k\to \Mod_k$ is equivalent to $\DD_2(A)$ not in 
$\Alg_2(\Mod_k)$
but in
$\Alg_1(\Mod_k)$.
But,the endomorphism of $k\in \Coh(k\otimes_Ak)$,
which we regard as the endomorphism algebra of the identity functor through $\Coh(k\otimes_Ak)\hookrightarrow \End_{\Ind(\textup{LCoh}(A))}(\Mod_k)$,
is equivalent to the endomorphism of the unit object $\mathsf{D}(\LMod_A)$
through the compatible fully faithful embedding $\Coh(k\otimes_Ak) \subset\LMod_{k\otimes_Ak}\simeq \mathsf{D}(\LMod_A)$.
By Lemma~\ref{dualcatalg},
the latter endomorphism $\etwo$-algebra
is equivalent to $\DD_2(A)$.
\QED

\begin{Remark}
\label{etwocaserem}
Let $A$ be an Artin $\etwo$-algebra.
Then the functor 
\[
\Phi^!:\Ind(\textup{LCoh}(k\otimes_Ak))_{<\infty}\to \End_{\Ind(\textup{LCoh}(A))}(\Mod_k)
\]
sending $M$ to the $!$-integral transform $q_*(p^!(-)\otimes^!M)$ is fully faithful.
Here $p^!:\Mod_k \to \Ind(\textup{LCoh}(k\otimes_Ak))$ is the $!$-pullback which is right adjoint 
to the forgetful functor along $k\to k\otimes_Ak$, and $q_*: \Ind(\textup{LCoh}(k\otimes_Ak))\to \Mod_k$
is the forgetful functor along $k\otimes_Ak\leftarrow k$.
 $\Ind(\textup{LCoh}(k\otimes_Ak))^{\le0}\to \End_{\Ind(\textup{LCoh}(A))}(\Mod_k)$
is fully faithful since the proof \cite[Theorem 4.0.5 (ii)]{BPN} works in this situation: the points required for the arguments are (i) $\LMod_B\otimes_k\LMod_C\simeq \LMod_{B\otimes_k C}$ for $\etwo$-algebras $B$, $C$ over $k$,
(ii) for a connective $\etwo$-algebra $B$ the colimit-preserving functor $\Ind(\textup{LCoh}(B))\to \textup{LMod}_B$ extending the inclusion
$\textup{LCoh}(B)\hookrightarrow \textup{LMod}_B$ is $t$-exact and commutes with pushforward functors, and it also induces
an equivalence $\Ind(\textup{LCoh}(B))^{\ge0}\to \textup{LMod}_B^{\ge0}$.
\end{Remark}

\begin{Construction}
As a consequence of Lemma~\ref{Indend}, the counit map $\Theta_k\circ E_k\to \textup{id}$ determines the canonical monoidal functor 
\[
\rho_A:\LMod_{\DD_2(A)}^\otimes \longrightarrow \End_{\Ind(\textup{LCoh}(A))}(\Mod_k)^\otimes.
\]
We will consider
the biduality map
\[
\delta_{\Ind(\textup{LCoh}(A))}:\Ind(\textup{LCoh}(A))\longrightarrow \mathsf{D}(\mathsf{D}(\Ind(\textup{LCoh}(A))))=\End_{\End_{\Ind(\textup{LCoh}(A))}(\Mod_k)}(\Mod_k)
\]
and the restriction along $\rho_A$
\[
U_A:\End_{\End_{\Ind(\textup{LCoh}(A))}(\Mod_k)}(\Mod_k)\to \End_{\LMod_{\DD_2(A)}}(\Mod_k).
\]
The composition gives rise to the monoidal functor
\[
\chi_A:\Ind(\textup{LCoh}(A))\longrightarrow \End_{\LMod_{\DD_2(A)}}(\Mod_k).
\]
\end{Construction}


\begin{Theorem}
\label{main1}
The morphism $\chi_A$ in $\Alg_1(\PR_k)$ is an equivalence.
\end{Theorem}

\begin{Remark}
We showed an equivalence $\Ind(\textup{LCoh}(A))\simeq \End_{\LMod_{\DD_2(A)}}(\Mod_k)$ in $\PR_k$ in \cite{HD},
but a monoidal equivalence is a much stronger result.
\end{Remark}

\section{$t$-structures and Duality results}

\label{DS}

In this section, we formulate Theorem~\ref{intro0} and Theorem~\ref{intro3}.
To this end, we take into account various $t$-structures and its left completeness and right completeness.
We begin by reviewing $t$-structures.

\subsection{Generalities on $t$-structures}

\label{tsec}

Unless otherwise stated, we adopt the cohomological index.
Let us recall the notion of $t$-structures on stable $\infty$-categories.
Let $\CCC$ be a stable $\infty$-category.
Recall that the homotopy category $\hhh(\CCC)$ has the structure of a triangulated category.
A $t$-structure on $\CCC$ is a pair of full subcategories $(\CCC^{\le0},\CCC^{\ge0})$ of $\CCC$
such that $(\hhh(\CCC^{\le0}), \hhh(\CCC^{\ge0}))$ is a $t$-structure on the triangulated category $\hhh(\CCC)$ (see e.g. \cite{HA}, \cite{KS}, \cite[Section 2.1]{IH}).
We shall write $\CCC^{\le n}$ and $\CCC^{\ge n}$ for $\CCC^{\le 0}[-n]$ and $\CCC^{\ge 0}[-n]$, respectively.

Suppose that $(\CCC^{\le0},\CCC^{\ge0})$ is a $t$-structure.
The inclusion $\CCC^{\le n}\hookrightarrow\CCC$ has a right adjoint 
$\tau^{\le n}:\CCC\longrightarrow \CCC^{\le n}$
and the inclusion $\CCC^{\ge n}\hookrightarrow\CCC$
has a left adjoint $\tau^{\ge n}:\CCC\longrightarrow \CCC^{\ge n}$.
There exists a canonical equivalence $\tau^{\ge 0}\circ \tau^{\le0}\simeq \tau^{\le 0}\circ \tau^{\ge0}$.
We often use the notation $(\CCC^{\le0}, \CCC^{\ge0})$ to denote a $t$-structure on $\CCC$.
Note that $\CCC^{\ge0}$ is uniquely determined by $\CCC^{\le0}$: the full subcategory $\CCC^{\ge0}$ is spanned by those objects $C'$ such that $\Map_{\CCC}(C,C')$ is contractible for any $C\in \CCC^{\le -1}$. Dually, $\CCC^{\le0}$ is uniquely determined by $\CCC^{\ge0}$.
By $(\CCC,\CCC^{\le0})$ we mean a stable $\infty$-category endowed with the $t$-structure determined by a full subcategory $\CCC^{\le0}$.

\begin{Definition}
Let $(\CCC^{\le 0}, \CCC^{\ge 0})$ be a $t$-structure on a stable $\infty$-category $\CCC$.
\begin{enumerate}
\item For any $n \in \ZZ$, the left adjoint functor $\tau^{\ge n}: \CCC \to \CCC^{\ge n}$ to the inclusion
$\CCC^{\ge n} \subset \CCC$ induces a morphism $C \to \tau^{\ge n}(C)$, determined by the unit map. For all $n \in \ZZ$, this yields the sequence
\[
C \to \cdots \to \tau^{\ge n-1}(C) \to \tau^{\ge n}(C) \to \tau^{\ge n+1}(C) \to \cdots
\]
in $\CCC$.
The $t$-structure is left complete if, for any $C \in \CCC$, the induced morphism $C \to \varprojlim_{n \in \ZZ} \tau^{\ge n}(C)$ is an equivalence.
\item The right adjoint functor $\tau^{\le n}: \CCC \to \CCC^{\le n}$ to the inclusion
$\CCC^{\le n} \subset \CCC$ induces a morphism $\tau^{\le n}(C) \to C$, determined by the counit map. For all $n \in \ZZ$, this yields the sequence
\[
\cdots \to \tau^{\le n-1}(C) \to \tau^{\le n}(C) \to \tau^{\le n+1}(C) \to \cdots \to C
\]
in $\CCC$.
The $t$-structure is right complete if, for any $C \in \CCC$, the induced morphism $\varinjlim_{n \in \ZZ} \tau^{\le n}(C) \to C$ is an equivalence.
\item The $t$-structure is left separated if $\bigcap_{n \in \ZZ} \CCC^{\le n} = 0$.
Dually, it is right separated if $\bigcap_{n \in \ZZ} \CCC^{\ge n} = 0$.
\item If $\CCC$ admits filtered colimits and the full subcategory $\CCC^{\ge 0}$ is closed under filtered colimits in $\CCC$, we say that the $t$-structure is compatible with filtered colimits.
\item Suppose that $\CCC$ is presentable. The $t$-structure $(\CCC^{\le 0}, \CCC^{\ge 0})$ is accessible if the full subcategory $\CCC^{\le 0}$ is presentable.
\item Suppose that $\CCC$ has a monoidal structure. If the tensor product functor $\CCC \times \CCC \to \CCC$ sends $\CCC^{\le 0} \times \CCC^{\le 0}$ to $\CCC^{\le 0}$ and the unit object lies in $\CCC^{\le 0}$, we say that the $t$-structure is compatible with the monoidal structure.
\end{enumerate}
\end{Definition}

\begin{Remark}
If $C\in \cap_{n\in \ZZ}\CCC^{\le n}$, then $\tau^{\ge n}C\simeq0$ for all $n\in \ZZ$.
It follows that if a $t$-structure $(\CCC^{\le0}, \CCC^{\ge0})$ is left complete,
then it is left separated.
Dually, if a $t$-structure is right complete,
then it is right separated.
\end{Remark}

\begin{Example}
\label{connectivetstr}
Let $A$ be an associative connective ring spectrum, and let $\LMod_A$ be the stable $\infty$-category of left $A$-module spectra.
Define $\LMod_A^{\le 0}$ as the full subcategory of $\LMod_A$ spanned by objects $M$ such that $\pi_n(M) = 0$ for all $n < 0$, and $\LMod_A^{\ge 0}$ as the full subcategory spanned by objects $M$ such that $\pi_n(M) = 0$ for all $n > 0$.
Then $(\LMod_A^{\le 0}, \LMod_A^{\ge 0})$ defines an accessible, left complete, and right complete $t$-structure (see \cite[7.1.1.13]{HA}).
The full subcategory $\LMod_A^{\le 0}$ is also characterized as the smallest full subcategory which contains $A$ and is closed under small colimits.
Moreover, the $t$-structure is compatible with filtered colimits, as $\LMod_A^{\ge 0} \subset \LMod_A$ is closed under filtered colimits.
This $t$-structure is compatible with the monoidal structure when $A$ is an $\etwo$-algebra.
The heart $\LMod_A^{\heartsuit}$ is equivalent to the abelian category of left $\pi_0(A)$-modules.
\end{Example}

\begin{Remark}
\label{unifiedtstr}
Example~\ref{connectivetstr} is generalized to the nonconnective case as follows.
For $A\in \Alg_1(\SP)$ we define $\LMod_A^{\le0}$ to be the smallest full subcategory
which contains $A$ and is closed under colimits and extensions. By \cite[1.4.4.11]{HA},
$\LMod_A^{\le0}$ determines an accessible $t$-structure on $\LMod_A$ such that $\LMod_A^{\ge0}$
consists of those objects $M$ such that $\pi_0(\Map_{\LMod_A}(A[n],M))=\pi_n(M)=H^{-n}(M)=0$ for $n>0$.
The full subcategory $\LMod_A^{\ge0}$ is closed under filtered colimits, and the $t$-structure is right complete.
We shall refer to this $t$-structure as the standard $t$-structure.
When $A$ is an $\etwo$-algebra, the $t$-structure is compatible with the monoidal structure on $\LMod_A$.
If $A$ is not connective, the resulting $t$-structure is generally not left complete.
\end{Remark}

\begin{Example}
\label{dualstand}
For $A\in \Alg_1(\SP)$ we define $\LMod_A^{\le0,ds}$ to be the fiber product
$\LMod_A\times_{\SP}\SP^{\le0}$ where $\LMod_A\to \SP$ is the (colimit-preserving) forgetful functor
and $\SP^{\le0}\to \SP$ is the inclusion which preserves small colimits.
This fiber product can also be regarded as a fiber product of presentable $\infty$-categories.
By definition, $\LMod_A^{\le0,ds}\subset \LMod_A$ is closed under extensions and small colimits
so that it determines an accessible $t$-structure (cf. \cite[1.4.4.11]{HA}).
We refer to this $t$-structure as the dual standard $t$-structure.
To distinguish this $t$-structure from one in Remark~\ref{unifiedtstr},
we denote it by $(\LMod_A^{\le0,ds},\LMod_A^{\ge0,ds})$.
We also define the $t$-structure $(\RMod_A^{\le0,ds}, \RMod_A^{\ge0,ds})$ on $\RMod_A$ in a similar way.
If $A$ is not connective, the resulting $t$-structure is generally not right complete (cf. Corollary~\ref{nonright}).
\end{Example}

\begin{Example}
\label{Indtstr}
Let $\CCC$ be a small stable $\infty$-category.
Let $(\CCC^{\le0}, \CCC^{\ge0})$ be a $t$-structure on $\CCC$.
Then $(\Ind(\CCC^{\le0}),\Ind(\CCC^{\ge0}))$ defines a $t$-structure on $\Ind(\CCC)$.
By the construction, the $t$-structure is compatible with filtered colimits.
If the $t$-structure on $\CCC$ is right bounded, then 
$(\Ind(\CCC^{\le0}),\Ind(\CCC^{\ge0}))$ is right complete (see \cite[VIII, Lemma 5.4.1]{DAG}).

For instance, suppose that $R$ is a noetherian connective $\eone$-algebra over a field $k$, 
that is, an object of $\Alg_1(\Mod_k^{\le0})$.
We let denote $\textup{LCoh}(R)$ the full subcategory of $\LMod_R$ spanned by those objects $M$
such that $H^n(M)=0$ for $|n|>>0$, and each $H^n(M)$ is a finitely generated $H^0(R)$-module.
This small stable $\infty$-category $\textup{LCoh}(R)$ inherits a bounded $t$-structure from that of $\LMod_R$ in Example~\ref{connectivetstr}.
Then $(\Ind(\textup{LCoh}(R)^{\le0}),\Ind(\textup{LCoh}(R)^{\ge0}))$ 
is a right complete accessible $t$-structure. This is generally not left complete.
\end{Example}

Let $\PRT$ be the $\infty$-category of presentable stable $\infty$-categories equipped with accessible $t$-structures.
An object of $\PRT$ is a pair $(\CCC,\CCC^{\le0})$ where $\CCC$ is a presentable stable $\infty$-category, and $\CCC^{\le0}$ is a  full subcategory
inducing a $t$-structure,
such that $\CCC^{\le0}$ itself is presentable.
A morphism $(\CCC, \CCC^{\le0})\to (\DDD,\DDD^{\le0})$ is a colimit-preserving functor $F:\CCC\to  \DDD$ which is right $t$-exact (i.e., $F(\CCC^{\le0})\subset \DDD^{\le0}$).
Let $\PRTT$ be the full subcategory of $\PRT$, spanned by presentable stable $\infty$-categories equipped with accessible right complete $t$-structures.
Let $\PRTTT$ be the full subcategory of $\PRT$, spanned by presentable stable $\infty$-categories equipped with accessible left complete and right complete $t$-structures.
There exists an adjoint pair
\[
\xymatrix{
I^+: \PRTT \ar@<0.5ex>[r] &  \PRT:\widehat{R}  \ar@<0.5ex>[l]  
}
\]
such that $I^+$ is the inclusion, and the right adjoint $\widehat{R}$ carries $(\CCC,\CCC^{\le0})$
to the right completion $\widehat{R}(\CCC,\CCC^{\le0})$ given by
$(\colim_{n\to \infty} \CCC^{\le n},\CCC^{\le0})$
where $\colim_{n\to \infty} \CCC^{\le n}$ is the colimit of the sequence of inclusions
\[
\CCC^{\le 0}\to \CCC^{\le 1}\to \CCC^{\le 2}\to \cdots\to \CCC^{\le n}\to \cdots
\]
in $\PR_{\mathbb{S}}$
(cf. \cite[5.5.3.4, 5.5.3.18]{HTT}). 

The left completion $\HL((\CCC,\CCC^{\le0}))$ of $(\CCC,\CCC^{\le0})$
is defined as
$\lim_{n\to -\infty}\CCC^{\ge n}$ in $\PR_{\mathbb{S}}$
of the sequence 
\[
\cdots \to \CCC^{\ge n-1}\to \CCC^{\ge n}\to \CCC^{\ge n+1} \to \cdots
\]
induced by truncations.
We usually write $\HL(\CCC)$ for $\HL((\CCC,\CCC^{\le0}))$.
If $(\CCC,\CCC^{\le0})\in \PRTT$, then the left completion $\HL(\CCC)$ is also accessible and right complete (see \cite[VIII, 4.6.13]{DAG}).
There exists a localization adjoint pair
\[
\xymatrix{
\HL :\PRTT \ar@<0.5ex>[r] &  \PRTTT:I^{\pm}  \ar@<0.5ex>[l]  
}
\]
such that $I^\pm$ is the inclusion, and $\HL$ is the left adjoint given by the left completion.

{\it Symmetric monoidal structures.}
Next, we recall the symmetric monoidal structure on $\PRT$ (see \cite[VIII, 4.6.1]{DAG}).
Let $(\mathcal{B},\mathcal{B}^{\le0})$ and $(\CCC,\CCC^{\le0})$ be presentable $\infty$-categories endowed with accessible $t$-structures. 
The tensor product of $(\mathcal{B},\mathcal{B}^{\le0})$ and $(\CCC,\CCC^{\le0})$
is defined as
$\mathcal{B}\otimes \CCC$ (in $\PR_{\mathbb{S}}$) 
together with the smallest full subcategory $(\mathcal{B}\otimes \CCC)^{\le0}$
which contains the essential image of $\mathcal{B}^{\le0}\times \CCC^{\le0}\to \mathcal{B}\otimes \CCC$
and is closed under small colimits and extensions.
By \cite[1.4.5.11]{HA}, $(\mathcal{B}\otimes \CCC)^{\le0}$ defines an accessible $t$-structure.
The universal property of $(\mathcal{B}\otimes \CCC, (\mathcal{B}\otimes \CCC)^{\le0})$ is described as follows: For any $(\DDD, \DDD^{\le0})\in \PRT$
the composition with the canonical functor $\mathcal{B}\times \CCC\to \mathcal{B}\otimes \CCC$ gives an equivalence
\[
\Fun_{\PRT}((\mathcal{B}\otimes \CCC, (\mathcal{B}\otimes \CCC)^{\le0}),(\DDD,\DDD^{\le0}))\to \Fun^{\le}(\mathcal{B}\times \CCC,\DDD)
\]
where $\Fun^{\le}(\mathcal{B}\times \CCC,\DDD)$ is the full subcategory of $\Fun(\mathcal{B}\times \CCC,\DDD)$,
which consists of those functors $F$ that preserve small colimits separately in each variable, such that the essential image of $\mathcal{B}^{\le0}\times \CCC^{\le0}$ 
is contained in $\DDD^{\le0}$.

If $(\mathcal{B},\mathcal{B}^{\le0})$ and $(\CCC,\CCC^{\le0})$ are right complete,
then the tensor product $(\mathcal{B}\otimes \CCC, (\mathcal{B}\otimes \CCC)^{\le0})$ in $\PRT$
is also right complete (see \cite[VIII, 4.6.11]{DAG}).
Thus, $\PRTT$ inherits a symmetric monoidal structure from that of $\PRT$.
The left adjoint $I^+$ is promoted
to a symmetric monoidal functor.
The $\infty$-category $\PRTTT$ also has a symmetric monoidal structure (see \cite[VIII, Proposition 4.6.15, Lemma 4.6.14]{DAG}).
Let $(\mathcal{B},\mathcal{B}^{\le0})$ and $(\CCC,\CCC^{\le0})$ be objects of $\PRTTT$.
Explicitly, the tensor product of $(\mathcal{B},\mathcal{B}^{\le0})$ and $(\CCC,\CCC^{\le0})$ in $\PRTTT$
is given by 
\[
\HL((\mathcal{B}\otimes\CCC,  (\mathcal{B}\otimes \CCC)^{\le0})).
\]
We usually write $\HL(\mathcal{B}\otimes\CCC)$ or $\mathcal{B}\widehat{\otimes}\CCC$ for it.
The functor $\HL:\PRTT\to \PRTTT$ is promoted to a symmetric monoidal functor. 
Since $\HL$ is a left adjoint functor, $\HL$ commutes with the relative tensor products.
The inclusion $I^{\pm}$
is not a symmetric monoidal functor but a lax symmetric monoidal functor.

\vspace{2mm}

For later use, we prove the following two lemmas.

\begin{Lemma}
\label{functstr}
Let $\mathcal{A}\in \Alg_1(\PRTT)$.
Let $\mathcal{C}\in \LMod_{\mathcal{A}}(\PRTT)$.
Let $\mathcal{D}$ be an object of $\LMod_{\mathcal{A}}(\PRTT)$.
Then $\Fun_{\mathcal{A}}(\CCC, \DDD)$ has an accessible $t$-structure such that $\Fun_{\mathcal{A}}(\CCC, \DDD)^{\le0}$
is the full subcategory of right $t$-exact functors.
\end{Lemma}

\begin{Remark}
The full subcategry $\Fun_{\mathcal{A}}(\CCC, \DDD)^{\le0}$
can naturally be regarded as a presentable $\infty$-category $\Fun_{\mathcal{A}^{\le0}}(\CCC^{\le0}, \DDD^{\le0})$
where $\Fun_{\mathcal{A}^{\le0}}(\CCC^{\le0}, \DDD^{\le0})$ is the presentable $\infty$-category equipped with
\[
\CCC^{\le0}\otimes\Fun_{\mathcal{A}^{\le0}}(\CCC^{\le0}, \DDD^{\le0})\to \DDD^{\le0}
\]
($\otimes$ is the tensor product in $\PR$)
such that for any $\PPP\in \PR$, the induced map
$\Map_{\PR}(\PPP, \Fun_{\mathcal{A}^{\le0}}(\CCC^{\le0}, \DDD^{\le0}))\to \Map_{\LMod_{\mathcal{A}^{\le0}}(\PR)}(\CCC^{\le0}\otimes\PPP, \DDD^{\le0})$
is an equivalence of spaces.
\end{Remark}

\Proof
Since $(\CCC^{\le0}, \CCC^{\ge0})$ is right complete, it follows from  \cite[VIII, Lemma 4.6.6]{DAG}
that the restriction to $\CCC^{\le0}\subset\CCC$ gives an equivalence
$\Fun_{\mathcal{A}}(\CCC,\DDD)\simeq \Fun_{\mathcal{A}^{\le0}}(\CCC^{\le0},\DDD)$.
Consequently, 
the presentable $\infty$-category
$\Fun_{\mathcal{A}^{\le0}}(\CCC^{\le0}, \DDD^{\le0})$ 
is a full subcategory $ \Fun_{\mathcal{A}^{\le0}}(\CCC^{\le0},\DDD)$ which is closed under small colimits.
Thus, if
$\Fun_{\mathcal{A}^{\le0}}(\CCC, \DDD)^{\le0}$ is defined to be $\Fun_{\mathcal{A}^{\le0}}(\CCC^{\le0}, \DDD^{\le0})$,
it determines an accessible $t$-structure.
\QED

\begin{Lemma}
\label{completionfunctor}
Let $\mathcal{A}\in \Alg_1(\PRTTT)$.
Let $\mathcal{C}\in \LMod_{\mathcal{A}}(\PRTT)$.
Let $\mathcal{D}$ be an object of $\LMod_{\mathcal{A}}(\PRTTT)$.
Let $\mathcal{C}\to \widehat{\CCC}:=\widehat{L}(\CCC)$ be the left completion of $\CCC$.
We abuse notation by writing $\Fun_{\mathcal{A}}(\widehat{\CCC},\DDD)$ and $\Fun_{\mathcal{A}}(\CCC,\DDD)$
for presentable $\infty$-categories of $\mathcal{A}$-linear functors where we discard $t$-structures.
We equip $\Fun_{\mathcal{A}}(\widehat{\CCC},\DDD)$  and $\Fun_{\mathcal{A}}(\CCC,\DDD)$
with accessible $t$-structures defined in Lemma~\ref{functstr}.
Then the composition induces an equivalence
$\widehat{R}(\Fun_{\mathcal{A}}(\widehat{\CCC},\DDD))\simeq \widehat{R}(\Fun_{\mathcal{A}}(\CCC,\DDD))$
(where $\widehat{R}$ is the right completion functor).
\end{Lemma}

\Proof
Since the $t$-structure on $\CCC$ is right  complete (it also holds for $\mathcal{A}$), then
by \cite[VIII, Lemma 4.6.6]{DAG}, the restriction along $\CCC^{\le0}\subset \CCC$ induces an equivalence
\begin{eqnarray*}
\Fun_{\mathcal{A}}(\CCC,\DDD) &\simeq&  \Fun_{\mathcal{A}^{\le0}}(\CCC^{\le0},\DDD).
\end{eqnarray*}
Moreover, by definition
\begin{eqnarray*}
\widehat{R}(\Fun_{\mathcal{A}^{\le0}}(\CCC^{\le0},\DDD))&\simeq& \colim_n\Fun_{\mathcal{A}^{\le0}}(\CCC^{\le0},\DDD^{\le n})  \\
\end{eqnarray*}
where the right-hand side is a colimit as a presentable stable $\infty$-category. 
Similarly, there is an equivalence
\[
\widehat{R}(\Fun_{\mathcal{A}}(\widehat{\CCC},\DDD))\simeq \colim_n\Fun_{\mathcal{A}^{\le0}}(\widehat{\CCC}^{\le0},\DDD^{\le n}).
\]
Since $\widehat{\CCC}$ is the left completion and $\DDD$ is left complete, 
we have an equivalence 
\[
\Fun_{\mathcal{A}^{\le0}}(\CCC^{\le0}, \DDD^{\le0})\simeq \Fun_{\mathcal{A}^{\le0}}(\widehat{\CCC}^{\le0},\DDD^{\le0}).
\]
It follows that
\[
\colim_n\Fun_{\mathcal{A}^{\le0}}(\CCC^{\le0},\DDD^{\le n}) \simeq \colim_n\Fun_{\mathcal{A}^{\le0}}(\widehat{\CCC}^{\le0},\DDD^{\le n}). 
\]
Consequently, we have $\widehat{R}(\Fun_{\mathcal{A}}(\widehat{\CCC},\DDD))\simeq \widehat{R}(\Fun_{\mathcal{A}}(\CCC,\DDD))$.
Unfolding constructions, we see that this is induced by the composition with $\CCC\to \widehat{\CCC}$.
\QED

\subsection{}

\begin{Notation}
Let $(\PRTT)_k$ denote $\Mod_{\Mod_k}(\PRTT)$ where we regard $\Mod_k$ as an object of $\CAlg(\PRTT)$ by taking $\Mod_k^{\le0}$
to be the full subcategory spanned by connective $k$-module spectra. 
We write $(\PRTTT)_k$ for $\Mod_{\Mod_k}(\PRTTT)$.
\end{Notation}

\begin{Notation}
\label{compD}
We consider $\mathcal{A}\in \Alg_1((\PRTT)_k)_{/\Mod_k}$ and $\mathsf{D}(\mathcal{A})=\End_{\mathcal{A}}(\Mod_k)$.
As we will observe in Lemma~\ref{functstr}, $\mathsf{D}(\mathcal{A})$ admits an accessible $t$-structure
determined by $\mathsf{D}(\mathcal{A})^{\le0}=\Fun_{\mathcal{A}^{\le0}}(\Mod_k^{\le0},\Mod_k^{\le0})$.
This $t$-structure is generally not right complete (see Corollary~\ref{nonright}): $\mathsf{D}$ does not generally work well. Thus, we define 
$\DDDD(\mathcal{A})$ to be the right completion of $\mathsf{D}(\mathcal{A})$.
\end{Notation}

\begin{Construction}
\label{diagramconst}
We have the diagram 
\[
\xymatrix{
\LMod_A \ar[r]^{\Upsilon_A} \ar[d]_{\delta_{\LMod_A}} &  \Ind(\textup{LCoh}(A)) \ar[d]_{\delta_{\Ind(\textup{LCoh}(A))}}   \ar[dr]^{\chi_A} &  \\
\mathsf{D}(\mathsf{D}(\LMod_A)) \ar[r]_{\mathsf{D}(\mathsf{D}(\Upsilon_A))} &  \mathsf{D}(\mathsf{D}(\Ind(\textup{LCoh}(A)))) \ar[r]^{U_A} &   \End_{\LMod_{\DD_2(A)}}(\Mod_k)
}
\]
where vertical functors $\delta_{(-)}$ are biduality functors, and $\Upsilon_A$ is the canonical monoidal functor given informally by the tensor with
the dualizing object $(-)\otimes_A\omega_A$ (see Appendix).
See Section~\ref{ICD}, Theorem~\ref{main1} for $U_A$, $\chi_A$. 
\end{Construction}

\begin{Theorem}
\label{main0}
There exists an equivalence
\[
\LMod_A \simeq \widehat{R}(\mathsf{D}(\mathsf{D}(\LMod_A)))
\]
in $\Alg_1(\PR_k)$.
\end{Theorem}

\begin{Remark}
Since the $t$-structure on $\mathsf{D}(\LMod_A)$ is right complete, we have $\widehat{R}(\mathsf{D}(\mathsf{D}(\LMod_A)))\simeq \DDDD(\DDDD(\LMod_A))$.
\end{Remark}

\begin{Remark}
We do not know if $\mathsf{D}(\mathsf{D}(\LMod_A))$ is right complete. 
\end{Remark}

\begin{Proposition}
\label{compositeequivalence}
The composite $U_A\circ \mathsf{D}(\mathsf{D}(\Upsilon_A))$ in Construction~\ref{diagramconst} is an equivalence after a right completion.
\end{Proposition}

The (symmetric) monoidal functor $\widehat{L}$
gives rise to 
the adjoint pair 
\[
\xymatrix{
\Alg_1(\PRTT) \ar@<0.3ex>[r]    &  \Alg_1(\PRTTT)  \ar@<0.3ex>[l] 
}
\]
such that the left adjoint carries $\mathcal{A}$ to $\widehat{L}(\mathcal{A})$
where $\mathcal{A}$ denotes an associative algebra object in $\PRTT$ (i.e., a presentably monoidal $\infty$-category with a compatible accessible right complete $t$-structure),
 and $\widehat{L}(\mathcal{A})$ denotes the left-completed one.
The right adjoint is the inclusion.
According to \cite[Lemma 3.10]{IH},
$\xi:\LMod_{\DD_2(A)}\to \End_{\LMod_A}(\Mod_k)$ defined as the counit map $\Theta_k\circ E_k(\End_{\LMod_A}(\Mod_k))\to \End_{\LMod_A}(\Mod_k)$ 
is a left completion. Namely, we have $\LMod_{\DD_2(A)}\to \widehat{L}(\LMod_{\DD_2(A)})\simeq \End_{\LMod_A}(\Mod_k)$
(in {\it loc. cit.} we write $\End_{A}^l(\Mod_k)$ for $\End_{\LMod_A}(\Mod_k)$).
We have 
\[
H:\End_{\End_{\LMod_A}(\Mod_k)}(\Mod_k) \to  \End_{\LMod_{\DD_2(A)}}(\Mod_k)
\]
induced by the restriction along $\xi$.

\vspace{3mm}

{\it Proof of Proposition~\ref{compositeequivalence}.}
We observe 
the diagram $(A)$:
\[
\xymatrix{
\Ind(\textup{LCoh}(k\otimes_Ak)) \ar[r]^{\Phi^!} \ar[d]^R &   \End_{\Ind(\textup{LCoh}(A))}(\Mod_k) \ar[d]^r \\
\LMod_{k\otimes_Ak} \ar[r]^{\Phi} &  \End_{\LMod_{A}}(\Mod_k)
}
\]
commutes up to homotopy,
where $R$ is the colimit-preserving functor extending $\Coh(k\otimes_Ak)\hookrightarrow \Mod_{k\otimes_Ak}$, and 
$\Phi^!$ and $\Phi$ are the $!$-integral transform $\Phi^!_K(-)=q_*(p^!(-)\otimes^! K)$ and the integral transforms $\Phi(K')=q_*(p^*(-)\otimes K')$.
This follows from the argument of \cite[Corollary 4.0.8]{BPN}.
Taking into account Lemma~\ref{Indend} and the commutativity of the diagram $(A)$,
we see that the composite $U_A\circ \mathsf{D}(\mathsf{D}(\Upsilon_A))$ can be identified with $H$.
It is enough to prove that $H$ is an equivalence after passing to the right completion $\widehat{R}$.
For ease of notation, we write $\mathcal{A}=\LMod_{\DD_2(A)}$, $\widehat{\mathcal{A}}=\widehat{L}(\LMod_{\DD_2(A)})\simeq \End_{\LMod_A}(\Mod_k)$
and $\mathcal{M}=\Mod_k$.
Note that $\mathcal{M}\in \PRTT$ lies in $\PRTTT$. 
First, there are canonical equivalences of $\infty$-categories 
\[
\End_{\widehat{\mathcal{A}}}(\mathcal{M},\mathcal{M})\simeq \Fun_{\mathcal{M}}(\mathcal{M}\otimes_{\widehat{\mathcal{A}}}\mathcal{M},\mathcal{M}) \ \ \ \ \textup{and}\ \ \ \ \End_{\mathcal{A}}(\mathcal{M},\mathcal{M})\simeq \Fun_{\mathcal{M}}(\mathcal{M}\otimes_{\mathcal{A}}\mathcal{M},\mathcal{M}).
\]
By Lemma~\ref{completionfunctor}, the compositions with the left completions
$\mathcal{M}\otimes_{\mathcal{A}}\mathcal{M}\to \widehat{L}(\mathcal{M}\otimes_{\mathcal{A}}\mathcal{M})\simeq \widehat{L}(\mathcal{M}\otimes_{\widehat{\mathcal{A}}}\mathcal{M})$ and $\mathcal{M}\otimes_{\widehat{\mathcal{A}}}\mathcal{M}\to \widehat{L}(\mathcal{M}\otimes_{\widehat{\mathcal{A}}}\mathcal{M})$
induce equivalences
\[
 \widehat{R}(\Fun_{\mathcal{M}}(\mathcal{M}\otimes_{\widehat{\mathcal{A}}}\mathcal{M},\mathcal{M})) \simeq  \widehat{R}(\Fun_{\mathcal{M}}(\widehat{L}(\mathcal{M}\otimes_{\widehat{\mathcal{A}}}\mathcal{M}),\mathcal{M}))  \simeq  \widehat{R}(\Fun_{\mathcal{M}}(\mathcal{M}\otimes_{\mathcal{A}}\mathcal{M},\mathcal{M})).
\]
By these equivalences, we see that the functor
$\widehat{R}(\End_{\widehat{\mathcal{A}}}(\mathcal{M}))\to  \widehat{R}(\End_{\mathcal{A}}(\mathcal{M}))$
 induced by the restriction along $\mathcal{A}\to \widehat{\mathcal{A}}$ is an equivalence.
\QED

Recall from Lemma~\ref{functstr} that $\Fun_{\LMod_{\DD_2(A)}}(\Mod_k,\Mod_k)$ has
 an accessible $t$-structure such that $\Fun_{\LMod_{\DD_2(A)}}(\Mod_k,\Mod_k)^{\le0}$
is defined as 
$\Fun_{\LMod_{\DD_2(A)}^{\le0}}(\Mod_k^{\le0},\Mod_k^{\le0})$.
The following gives a simple description of $\Fun_{\LMod_{\DD_2(A)}}(\Mod_k,\Mod_k)^{\le0}$.

\begin{Proposition}
\label{specialfunctstr}
There is an equivalence
$\Fun_{\LMod_{\DD_2(A)}}(\Mod_k,\Mod_k)\simeq \RMod_{\DD_1(A)}$.
In addition, this equivalence identifies the full subcategory
$\Fun_{\LMod_{\DD_2(A)}}(\Mod_k,\Mod_k)^{\le0}$ with $\RMod_{\DD_1(A)}^{\le0,ds}$ (cf. Example~\ref{dualstand}).
This is the smallest full subcategory $\langle k \rangle\subset \RMod_{\DD_1(A)}$ 
which contains the $\DD_1(A)$-module $k$ and is closed under small colimits.
\end{Proposition}

\begin{Remark}
\label{specialrem}
The object $k\in \RMod_{\DD_1(A)}$ corresponds to the identity functor $\Mod_k\to \Mod_k$.
\end{Remark}

\begin{Remark}
Since $\DD_2(A)$ is an $\etwo$-algebra there exists an equivalence $\DD_2(A)\simeq \DD_2(A)^{op}$
which induces $\LMod_{\DD_2(A)}\simeq \RMod_{\DD_2(A)}$.
Thus, one may replace $\RMod_{\DD_1(A)}$ with $\LMod_{\DD_2(A)}$ in the statement.
\end{Remark}

\begin{Corollary}
\label{tstromega}
Through the equivalence $\Ind(\textup{LCoh}(A))\simeq \End_{\LMod_{\DD_2(A)}}(\Mod_k)$ in Theorem~\ref{main1}, 
the full subcategory
$\End_{\LMod_{\DD_2(A)}}(\Mod_k)^{\le0}$ corresponds to the smallest full subcategory of 
$\Ind(\textup{LCoh}(A))$ which contains the unit object $\omega_A$ and is closed under small colimits and extensions.
\end{Corollary}

\begin{Definition}
\label{dualtind}
Let $\Ind(\textup{LCoh}(A))^{\le0,ds}$ denote the smallest full subcategory of $\Ind(\textup{LCoh}(A))$ which contains
$\omega_A$ and is closed under small colimits and extensions.
 By Corollary~\ref{tstromega}, it determines
an accessible $t$-structure.
We refer to it as the dual standard $t$-structure.
We often write $\Ind(\textup{LCoh}(A))^{ds}$ for $(\Ind(\textup{LCoh}(A)),\Ind(\textup{LCoh}(A))^{\le0,ds})$.
We note that this $t$-structure generally differs from that in \cite{Ind} (see Example~\ref{Indtstr}).
We refer to the latter $t$-structure determined by the connective part $\Ind(\textup{LCoh}(A)^{\le0})$ as the standard $t$-structure.  
\end{Definition}

\begin{Lemma}
\label{standardcompletion}
The functor
$\rho:\LMod_{\DD_1(A)}\to \RMod_A$ defined by $k\otimes_{\DD_1(A)}(-)$
is the left completion $\LMod_{\DD_1(A)}\to \widehat{L}(\LMod_{\DD_1(A)})$. Here
the $t$-structure on $\LMod_{\DD_1(A)}$ is the standard $t$-structure given in Remark~\ref{unifiedtstr}.
\end{Lemma}

\Proof
Since the right adjoint of $\rho$ is given by $\Hom_{\RMod_A}(k,-):\RMod_k\to \LMod_{\DD_1(A)}$, $\rho$ induces the equivalence
$\textup{LPerf}_{\DD_1(A)}\simeq \textup{RCoh}(A)$ where $\textup{LPerf}_{\DD_1(A)}$ is the smallest stable subcategory
which contains $\DD_1(A)$ and is closed under retracts (cf. Lemma~\ref{koszulfunc}).
Thus, $\LMod_{\DD_1(A)}\simeq \Ind(\textup{RCoh}(A))$, and through this equivalence
$\rho$ can be identified with the colimit-preserving functor $\Ind(\textup{RCoh}(A))\to \RMod_A$ that extends $\textup{RCoh}(A)\hookrightarrow \RMod_A$.
The argument in \cite[Proposition 1.2.6]{Ind} or \cite[Lemma 3.10]{IH} shows that $\Ind(\textup{RCoh}(A))\to \RMod_A$ is a left completion.
\QED

{\it Proof of Proposition~\ref{specialfunctstr}.}
There exist equivalences
\begin{eqnarray*}
\Fun_{\LMod_{\DD_2(A)}}(\Mod_k,\Mod_k) &\simeq& \Fun_k(\Mod_k\otimes_{\LMod_{\DD_2(A)}}\Mod_k,\Mod_k) \\
&\simeq& \Fun_k(\LMod_{k\otimes_{\DD_2(A)}k},\Mod_k) \\
&\simeq& \Fun_k(\LMod_{\DD_1(A)},\Mod_k) \\
&\simeq& \Fun_{\Mod_k^{\le0}}(\LMod_{\DD_1(A)}^{\le0,st},\Mod_k)
\end{eqnarray*}
in $\PR_k$ where the last equivalence is induced by the inclusion $\LMod_{\DD_1(A)}^{\le0,st}\subset \LMod_{\DD_2(A)}$
(see \cite[VIII, Lemma 4.6.6]{DAG}). 
Here the superscript $st$ indicates the $t$-structure in Remark~\ref{unifiedtstr} and stands for ``standard", which differs from
the dual standard $t$-structure (Example~\ref{dualstand}) in the coconnective case.
The third equivalence is induced by the equivalence $k\otimes_{\DD_2(A)}k\simeq \DD_1(A)$ (see \cite[Lemma 3.12]{HD}).
We observe that $\Fun_{\LMod_{\DD_2(A)}}(\Mod_k,\Mod_k)^{\le0}$ 
corresponds to $\Fun_{\Mod_k^{\le0}}(\LMod_{\DD_1(A)}^{\le0,st},\Mod_k^{\le0})$
because
$\Fun_{\LMod_{\DD_2(A)}^{\le0}}(\Mod_k^{\le0},\Mod_k^{\le0})\simeq \Fun_{\Mod_k^{\le0}}(\LMod_{\DD_1(A)}^{\le0,st},\Mod_k^{\le0})$.
Since $\RMod_{\DD_1(A)}$ is a dual of $\LMod_{\DD_1(A)}$ in $\PR_k$,
there exists an equivalence
$\Fun_k(\LMod_{\DD_1(A)},\Mod_k)\simeq \RMod_{\DD_1(A)}$.
This equivalence carries $F:\LMod_{\DD_1(A)}\to \Mod_k$
to $F(\DD_1(A))\in  \Mod_k$ endowed with the canonically defined right $\DD_1(A)$-module structure.
Thus, $\Fun_{\Mod_k^{\le0}}(\LMod_{\DD_1(A)}^{\le0,st},\Mod_k^{\le0})$ corresponds to
$\RMod_{\DD_1(A)}^{\le0, ds}$ (see Example~\ref{dualstand}).
Consequently, we have $\Fun_{\LMod_{\DD_2(A)}}(\Mod_k,\Mod_k)^{\le0}\simeq \RMod_{\DD_1(A)}^{\le0, ds}$.
According to Lemma~\ref{standardcompletion},
$\rho:\LMod_{\DD_1(A)}\to \RMod_A$
is a left completion so that
it induces an equivalence
\[
\Fun_{\Mod_k^{\le0}}(\LMod_{\DD_1(A)}^{\le0,st},\Mod_k^{\le0})\simeq \Fun_{\Mod_k^{\le0}}(\RMod_A^{\le0},\Mod_k^{\le0}).
\]
Consider the canonical equivalences
\begin{eqnarray*}
 \Fun_{\Mod_k^{\le0}}(\RMod_A^{\le0},\Mod_k)&\simeq&  \Fun_{\Mod_k}(\RMod_A,\Mod_k) \\
&\simeq& \LMod_A.
\end{eqnarray*}
The composite carries $F:\RMod_A^{\le0}\to \Mod_k$ to $F(A) \in \Mod_k$ endowed with the induced left $A$-module structure.
It follows that $\Fun_{\Mod_k^{\le0}}(\RMod_A^{\le0},\Mod_k^{\le0})$ corresponds to the full subcategory 
$\LMod_A^{\le0}$ through the equivalence.
Therefore, we have
\[
\RMod_{\DD_1(A)}^{\le0, ds}\simeq \Fun_{\LMod_{\DD_2(A)}}(\Mod_k,\Mod_k)^{\le0}\simeq \Fun_{\Mod_k^{\le0}}(\RMod_A^{\le0},\Mod_k^{\le0})\simeq \LMod_A^{\le0}
\]
which carries $k\in \RMod_{\DD_1(A)}^{\le0, ds}$ to $A\in \LMod_A^{\le0}$.
According to \cite[7.1.1.13]{HA}, $\LMod_A^{\le0}$ is
the smallest full subcategory which contains $A$ and is closed under small
colimits.
Thus, the presentable $\infty$-category $\RMod_{\DD_1(A)}^{\le0, ds}$ itself is the smallest
full subcategory which contains $k$ and closed under small colimits.
This completes the proof.
\QED

\begin{Corollary}
\label{nonright}
When $A$ is the trivial square zero extension $k\oplus k\cdot \epsilon=k[\epsilon]$ by $k\cdot\epsilon=k[1]$, the $t$-structure on 
$\Fun_{\LMod_{\DD_2(A)}}(\Mod_k,\Mod_k)$ is not right separated. In particular, it is not right complete.
\end{Corollary}

\Proof
If $A$ is $k\oplus k\cdot \epsilon=k[\epsilon]$, $\DD_1(A)$ is equivalent to a free associative algebra generated by $k[-2]$, which we denote by
$k[t]$ (cf. \cite[X, 4.4.5, 4.5.6]{DAG}).
Here, the cohomological degree of a generator $t$ is two.
In view of Proposition~\ref{specialfunctstr},
the $t$-structure on $\Fun_{\LMod_{\DD_2(A)}}(\Mod_k,\Mod_k)$
can be identified with
$(\RMod_{k[t]}^{\le0,ds}, \RMod_{k[t]}^{\ge0,ds})$.
Note also that $\RMod_{k[t]}^{\ge1,ds}$ is spanned those objects $M$ such that $\Map_{\RMod_{k[t]}}(k,M)$ is a contractible space.
The $k[t]$-module $k$ is a cofiber of the map $k[t][-2]\to k[t]$ given by the multiplication by $t$.
It follows that $\RMod_{k[t]}^{\ge1,ds}$ is spanned those objects $M$ such that
the fiber of the map $M\to  M[2]$ induced by the multiplication by $t$
is zero after taking the connective part $\tau^{\le0}:\Mod_k\to \Mod_k^{\le0}$.
Let $k[t,t^{-1}]$ be the $k[t]$-module obtained from $k[t]$ by inverting $t$
(more precisely, it is the colimit of $k[t]\hookrightarrow \frac{1}{t}k[t]\hookrightarrow  \frac{1}{t^2}k[t]\hookrightarrow \cdots$).
Then the multiplication by $t$ induces an equivalence $k[t,t^{-1}]\to k[t,t^{-1}][2]$ so that
$k[t,t^{-1}]$ belongs to $\bigcap_{n\in \ZZ}\RMod_{k[t]}^{\ge n,ds}$.
Thus, $\bigcap_{n\in \ZZ}\RMod_{k[t]}^{\ge n,ds}\neq \{0\}$ so that the $t$-structure on $\Fun_{\LMod_{\DD_2(A)}}(\Mod_k,\Mod_k)$
is not right separated. In particular, it is not right complete.
\QED

\section{Proofs of Duality results}

\label{proofsect}

In this section we give the proofs of Theorem~\ref{main1}, Theorem~\ref{main0}, and Theorem~\ref{intro3}.

\begin{Lemma}
\label{dualkernel}
Let $B,C$ be objects of $\Alg_1(\Mod_k)$.
Let $M$ be a right $B\otimes_kC$-module, that is, an object of $\RMod_{B\otimes_kC}$.
Let $f:\LMod_B\to \RMod_C$ be the functor given informally by $P\mapsto M\otimes_BP$,
where $M$ is regarded as a $C^{op}$-$B$-bimodule.
Then $f^\vee:\RMod_C^\vee\simeq \LMod_C\to \RMod_B \simeq \LMod_B^\vee$ induced by $f$
is determined by the integral kernel $M$. That is, $f^\vee$ is given informally by
$Q\mapsto M\otimes_CQ$ where $M$ is regarded as a $B^{op}$-$C$-bimodule.
\end{Lemma}

\Proof
Note that both $\LMod_B$ and $\RMod_C$ are dualizable in $\PR_k$.
Let $g:\Mod_k\to \LMod_B^\vee\otimes_k\RMod_C$
be the functor corresponding to $f$.
We also note that $f^\vee$ corresponds to $g:\Mod_k\to \LMod_B^\vee\otimes_k\RMod_C\simeq \LMod_B^{\vee}\otimes_k\RMod_C^{\vee\vee}$.
Through the canonical equivalences
\[
 \LMod_B^{\vee}\otimes_k\RMod_C\simeq  \RMod_B\otimes_k\RMod_C\simeq \RMod_{B\otimes_kC},
\]
$g$ can be identified with $\Mod_k\to \RMod_{B\otimes_kC}$, which is given by $T\mapsto T\otimes_kM$ where $M\in \RMod_{B\otimes_kC}$.
From these correspondences, we see that $f^\vee$
is given informally by
$Q\mapsto M\otimes_CQ$.
\QED

{\it Proof of Theorem~\ref{main1}.}
Let $\mu_A:\LMod_{\DD_2(A)}^\otimes\to \mathsf{D}(\LMod_A)$ be the algebraization functor (cf. Definition~\ref{algdef}, Lemma~\ref{dualcatalg}).
This is also equivalent to the morphism determined by the monoidal functor $P:\LMod_A\otimes_k\LMod_{\DD_2(A)}\simeq \LMod_{A\otimes_k\DD_2(A)}\to \Mod_k$
induced by base change along the universal pairing $A\otimes_k\DD_2(A)\to k$.
By Lemma~\ref{Indend} and the proof, we also have the algebraization functor $\rho_A:\LMod_{\DD_2(A)}^\otimes\to \mathsf{D}(\Ind(\textup{LCoh}(A)))$
such that $\mu_A$ factors as 
\[
\LMod_{\DD_2(A)}^\otimes \stackrel{\rho_A}{\longrightarrow} \mathsf{D}(\Ind(\textup{LCoh}(A))) \stackrel{\mathsf{D}(\Upsilon_A)}{\longrightarrow} \mathsf{D}(\LMod_A)
\]
where $\mathsf{D}(\Upsilon_A)$ is induced by $\Upsilon_A$.
Consider the adjoint pair
\[
\xymatrix{
\mathsf{D}: \Alg_1^+(\PR_k) \ar@<0.5ex>[r] & \Alg_1^+(\PR_k)^{op} :\mathsf{D} \ar@<0.5ex>[l]
}
\]
Using this adjoint pair, we have the commutative diagram 
\[
\xymatrix{
\LMod_A \ar[rr]^{\xi_A} \ar[rd]_{\Upsilon_A} & &  \mathsf{D}(\LMod_{\DD_2(A)}) \\
  & \Ind(\textup{LCoh}(A)) \ar[ru]_{\chi_A}  &  
}
\]
in $\Alg_1(\PR_k)$, where $\xi_A$ and $\chi_A$ correspond to $\mu_A$ and $\rho_A$, respectively.
Note that $\chi_A$ is determined by the pairing 
\[
\Ind(\textup{LCoh}(A))\otimes_k \LMod_{\DD_2(A)}^\otimes \stackrel{\textup{id}\otimes \rho_A}{\longrightarrow} \Ind(\textup{LCoh}(A))\otimes_k\mathsf{D}(\Ind(\textup{LCoh}(A))) \stackrel{\textup{universal}}{\longrightarrow}  \Mod_k.
\]
Likewise, $\xi_A$ is determined by $P$. Since $A$ is Artin, the biduality morphism $A\to \DD_2(\DD_2(A))$ is an equivalence.
It follows that $\xi_A$ is also an algebraization functor.

Recall from Proposition~\ref{specialfunctstr} that there exists an equivalence
$\mathsf{D}(\LMod_{\DD_2(A)})\simeq \RMod_{\DD_1(A)}$
in $\PR_k$.
Observe that $\xi_A$ carries $A$ to the identity functor $\Mod_k\to \Mod_k$ in $\mathsf{D}(\LMod_{\DD_2(A)})$, and the above equivalences send
the identity functor to $k\in \RMod_{\DD_1(A)}$ (coming from the augmentation $\DD_1(A)\to k$).
Passing to the endomorphism algebras,
$\xi_A$ induces an equivalence
$q:\End_{\LMod_A}(A)=A^{op} \to \End_{\RMod_{\DD_1(A)}}(k)=\End_{\LMod_{\DD_1(A)^{op}}}(k)\simeq \DD_1(\DD_1(A)^{op})=\DD_1(\DD_1(A))^{op}$.
Consequently, after discarding the monoidal structures, the underlying functor of $\xi_A$ can be identified with $k\otimes_{A}(-):\LMod_A\to \RMod_{\DD_1(A)}$,
where $k$ is the integral kernel which has the $\DD_1(A)^{op}$-$A$-bimodule structure arising from the universal pairing
$A\otimes_k\DD_1(A)\stackrel{q^{op}\otimes\textup{id}}{\to} \DD_1(\DD_1(A))\otimes_k\DD_1(A)\simeq A\otimes_k\DD_1(A)\to k$.

To prove our assertion, it is enough to show that the composite
$\Ind(\textup{LCoh}(A))\to \mathsf{D}(\LMod_{\DD_2(A)})\simeq\RMod_{\DD_1(A)}$ is an equivalence of plain $\infty$-categories.
Stable $\infty$-categories $\LMod_A$, $\Ind(\textup{LCoh}(A))$ and $\RMod_{\DD_1(A)}$
are compactly generated so that they are dualizable in $\PR_k$.
Duals of $\LMod_A$, $\Ind(\textup{LCoh}(A))$ and $\RMod_{\DD_1(A)}$
are $\RMod_A$, $\Ind(\textup{RCoh}(A))$ and $\LMod_{\DD_1(A)}$, respectively.
By \cite[4.8.4.1]{HA}, $\LMod_A^\vee\simeq \RMod_A$ and $\RMod_A^\vee\simeq \LMod_A$.
By Lemma~\ref{dualind} (2), we have $\Ind(\textup{LCoh}(A))^\vee\simeq \Ind(\textup{RCoh}(A))$.
Here the superscript $\vee$ indicates the dual in $\PR_k$.
We consider the diagram
\[
\xymatrix{
\LMod_{\DD_1(A)}\ar[rd]_{\chi_A^\vee} \ar[rr]^{\xi_A^\vee} &  & \RMod_A \\
 &  \Ind(\textup{RCoh}(A)) \ar[ru]_{\Psi_A} &
}
\]
obtained from 
$\LMod_A\to \Ind(\textup{LCoh}(A))\to \RMod_{\DD_1(A)}$
by taking duals.
By Construction~\ref{enconst},
$\Psi_A=\Upsilon_A^\vee$.
It suffices to prove that $\chi_A^\vee$ is an equivalence of $\infty$-categories.
To this end, we take into account $t$-structures.
Both $\RMod_A$ and $\LMod_{\DD_1(A)}$ are endowed with the standard $t$-structures (see Example~\ref{connectivetstr} and Remark~\ref{unifiedtstr}).
The stable $\infty$-category $\Ind(\textup{RCoh}(A))$ is equipped with the $t$-structure in Example~\ref{Indtstr}.
If we write 
$\pi:\RMod_A\to \Mod_k$ for the forgetful functor, $\pi^{-1}(\Mod_k^{\le0})=\RMod_A^{\le0}$ and $\pi^{-1}(\Mod_k^{\ge0})=\RMod_A^{\ge0}$.
Since the composite $\LMod_{\DD_1(A)}\to \RMod_A\stackrel{\pi}{\to} \Mod_k$  is
given by $k\otimes_{\DD_1(A)}(-)$, it follows that $\xi_A^\vee$ is $t$-exact.
By Lemma~\ref{dualkernel},
the functor $\xi_A^\vee$ is given by $k\otimes_{\DD_1(A)}(-)$ where $k$ is integral kernel having the $A^{op}$-$\DD_1(A)$-bimodule structure 
determined by the universal pairing.
Thus, we see that $\xi_A^\vee$ is $\rho$ in Lemma~\ref{standardcompletion}.
Lemma~\ref{standardcompletion} implies that $\xi_A^\vee$ is also a left completion functor.
In particular, $\xi_A^\vee$ induces an equivalence $\LMod_{\DD_1(A)}^{\ge0}\simeq  \RMod_A ^{\ge0}$.
Since $\Psi_A$ is a left completion functor (see \cite[Proposition 1.2.4]{Ind} for the commutative case, the $\etwo$-case is similar, cf. \cite[Lemma 3.11]{IH}),
$\Psi_A$ induces an equivalence $\Ind(\textup{RCoh}(A))^{\ge0}\simeq  \RMod_A ^{\ge0}$.
By Lemma~\ref{texactchidual} below, $\chi_A^\vee$
is also $t$-exact.
Thus, $\chi_A^\vee$ induces $\LMod_{\DD_1(A)}^{\ge0}\simeq \Ind(\textup{RCoh}(A))^{\ge0}$.
Note that there exist equivalences
$\Ind(\textup{RCoh}(A))^{\ge0}\simeq \Ind(\textup{RCoh}(A)^{\ge0})$, and $\LMod_{\DD_1(A)}^{\ge0}\simeq \Ind(\textup{RCoh}(A))^{\ge0}\simeq \RMod_A^{\ge0}$
(see the proof of Lemma~\ref{standardcompletion} for the latter equivalence).
Passing to compact objects, the restriction $\chi_A^{\vee}|_{\LMod_{\DD_1(A)}^{\ge0}}$ induces an equivalence $(\LMod_{\DD_1(A)}^{\ge0})^{\omega}\simeq \textup{RCoh}(A)^{\ge0}$. 
By the equivalence $\LMod_{\DD_1(A)}\simeq \Ind(\textup{RCoh}(A))$ (see Lemma~\ref{koszulfunc}),
any compact object of $\LMod_{\DD_1(A)}$ lies in $(\LMod_{\DD_1(A)}^{\ge0})^{\omega}$ up to a finite shift.
Similarly, any compact object of $\Ind(\textup{RCoh}(A))$ lies in $\textup{RCoh}(A)^{\ge0}$ up to a finite shift.
Consequently, we see that $\chi_A^{\vee}$ induces an equivalence  $\LMod_{\DD_1(A)}^{\omega}\simeq \textup{RCoh}(A)$
between the full subcategories of compact objects.
Since $\chi_A^\vee$ preserves small colimits, $\chi_A^{\vee}$ is an equivalence $\LMod_{\DD_1(A)}\simeq \Ind(\textup{RCoh}(A))$.
Finally, we give another argument for the final portion:
Both $\LMod_{\DD_1(A)}$ and $\Ind(\textup{RCoh}(A))$
are equipped with regular $t$-structures in the sense of \cite[Definition 4.3.1]{PR},
and they admit the equivalent left completions.
It follows from the uniqueness of regularization \cite[Proposition 4.3.2, Theorem 4.3.6]{PR}
that $\chi_A^\vee$ is an equivalence in $\PRT$.
\QED

\begin{Lemma}
\label{texactchidual}
Let $A$ be an Artin $\etwo$-algebra.
Let $\chi_A^\vee:\LMod_{\DD_1(A)}\simeq \RMod_{\DD_1(A)}^\vee\to \Ind(\textup{LCoh}(A))^\vee\simeq \Ind(\textup{RCoh}(A))$
denote 
the dual of $\chi_A:\Ind(\textup{LCoh}(A))\to \mathsf{D}(\LMod_{\DD_2(A)})\simeq \RMod_{\DD_1(A)}$
where $\mathsf{D}(\LMod_{\DD_2(A)})\simeq \RMod_{\DD_1(A)}$ is the equivalence in Proposition~\ref{specialfunctstr},
which sends the identity functor to $k\in \RMod_{\DD_1(A)}$.
Then $\chi_A^\vee$ carries $\DD_1(A)$ to $k\in \textup{RCoh}(A)\subset \Ind(\textup{RCoh}(A))$.
Moreover, $\chi_A^\vee$ is $t$-exact when 
$\LMod_{\DD_1(A)}$ is equipped with the standard $t$-structures (see Remark~\ref{unifiedtstr}), and
$\Ind(\textup{RCoh}(A))$ is equipped with the (standard) $t$-structure in Example~\ref{Indtstr}.
\end{Lemma}

\Proof
Note first that 
$\chi_A:\Ind(\textup{LCoh}(A))\to \mathsf{D}(\LMod_{\DD_2(A)})$
is determined by the pairing 
$\Ind(\textup{LCoh}(A))\otimes \LMod_{\DD_2(A)} \to \Mod_k$, which is the composite of $\textup{id}\otimes \rho_A$ and  the universal pairing (cf. the above proof).
Thus, if $u:\RMod_{\DD_1(A)}\to \Mod_k$ denotes the forgetful functor, the composite functor
$u\circ \chi_A:\Ind(\textup{LCoh}(A))\to \mathsf{D}(\LMod_{\DD_2(A)})\simeq \RMod_{\DD_1(A)}\to \mathsf{D}(\Mod_{k}) \simeq\Mod_k$
is equivalent to the $!$-pullback functor $\Ind(\textup{LCoh}(A))\to \Mod_k$ along the augmentation $A\to k$ (geometrically
corresponding to $\Spec k\to \Spec A$).
By Lemma~\ref{dualind} (3), this $!$-pullback functor is the
dual of the pushforward functor $f:\Mod_k\to \Ind(\textup{RCoh}(A))$;
$f$ is a unique colimit-preserving functor 
extending the restriction functor $\Perf_k\to \textup{RCoh}(A)\subset \Ind(\textup{RCoh}(A))$ 
along the augmentation $A\to k$.
Since both $\Ind(\textup{LCoh}(A))$ and $\RMod_{\DD_1(A)}$ are dualizable in $\PR_k$,
we see that $\chi_A^\vee\circ u^\vee:\Mod_k\to \LMod_{\DD_1(A)}\to \Ind(\textup{RCoh}(A))$
is equivalent to the pushforward functor $f$.
By Lemma~\ref{dualkernel}, the dual $u^{\vee}:\Mod_k\to \LMod_{\DD_1(A)}$ of the forgetful functor
is given informally by $T\mapsto \DD_1(A)\otimes_kT$. Namely, $u^\vee$ is a free functor.
In particular, $k\in \Mod_k$ maps to the free module $\DD_1(A)\in \LMod_{\DD_1(A)}$.
It follows that $\chi_A^\vee$ sends $\DD_1(A)$ to $k\in \Ind(\textup{RCoh}(A))$.
Next, we prove the latter assertion.
By the equivalence $h:\Ind(\textup{RCoh}(A))\simeq \LMod_{\DD_1(A)}$ in $\PRT$ (cf. Lemma~\ref{koszulfunc}; both are endowed with standard $t$-structures),
the $t$-structure on $\LMod_{\DD_1(A)}$ is right complete and is compatible with filtered colimits.
In addition, its heart $\LMod_{\DD_1(A)}^\heartsuit$ is equivalent to the abelian category
of (ordinary) $H^0(A)$-modules. Since $H^0(A)$ is an Artin local $k$-algebra with residue field $k$
(in particular, $\LMod_{\DD_1(A)}^\heartsuit$
is the Ind-category of the category of finitely generated $H^0(A)$-modules), to prove that $\chi_A^\vee$ is $t$-exact,
it is enough to show that $\chi_A^\vee$ carries $\DD_1(A)= h(k)\in \LMod_{\DD_1(A)}$ 
to an object of the heart in $\Ind(\textup{RCoh}(A))$.
This follows from the first assertion.
\QED

\begin{Remark}
Through the equivalence $\chi_A:\Ind(\textup{LCoh}(A))\simeq \mathsf{D}(\textup{LMod}_{\DD_2(A)})$,
the canonically defined $t$-structure on $\mathsf{D}(\textup{LMod}_{\DD_2(A)})$ (see Notation~\ref{compD} and Lemma~\ref{functstr})
corresponds to the dual $t$-structure on $\Ind(\textup{LCoh}(A))$ (cf. Definition~\ref{dualtind}, Proposition~\ref{specialfunctstr}, Remark~\ref{specialrem}).
In particular, $\chi_A$ is promoted to an equivalence in $\Alg_1((\PRT)_k)$.
\end{Remark}

Next, we prove 
Theorem~\ref{main0}.

{\it Proof of Theorem~\ref{main0}.}
The monoidal functor $\Upsilon_A:\LMod_A \to \Ind(\textup{LCoh}(A))$ is a fully faithful functor (see Proposition~\ref{pullbackdual}),
which induces an equivalence $\LMod_A^{\le0}\simeq \Ind(\textup{LCoh}(A))^{\le0,ds}$.
Moreover, the $t$-structure on $\LMod_A$ is right complete (see Example~\ref{connectivetstr}).
Thus, $\Upsilon_A$ induces $\LMod_A\simeq \widehat{R}(\Ind(\textup{LCoh}(A)))$.
By Theorem~\ref{main1}, there exists an equivalence 
\[
\widehat{R}(\Ind(\textup{LCoh}(A)))\simeq\widehat{R}(\End_{\LMod_{\DD_2(A)}}(\Mod_k)).
\]
By Proposition~\ref{compositeequivalence}, there exists an equivalence
\[
\widehat{R}(\mathsf{D}(\mathsf{D}(\LMod_A))\simeq \widehat{R}(\End_{\LMod_{\DD_2(A)}}(\Mod_k)).
\]
Finally, passing to the right completion to the diagram in Construction~\ref{diagramconst},
we obtain the diagram
\[
\xymatrix{
\LMod_A \ar[r]^{\simeq} \ar[d] &  \widehat{R}(\Ind(\textup{LCoh}(A))) \ar[d]   \ar[dr]^{\simeq} &  \\
\DDDD(\mathsf{D}(\LMod_A)) \ar[r] &  \DDDD(\mathsf{D}(\Ind(\textup{LCoh}(A)))) \ar[r] &  \widehat{R}(\End_{\LMod_{\DD_2(A)}}(\Mod_k)).
}
\]
We then see that the vertical arrow $\LMod_A\to \DDDD(\DDDD(\LMod_A)))\simeq \widehat{R}(\mathsf{D}(\mathsf{D}(\LMod_A)))$ is an equivalence.
\QED

\begin{Remark}
Through the equivalence $\Ind(\textup{LCoh}(A))\simeq \RMod_{\DD_1(A)}$ in Lemma~\ref{koszulfunc} (1),
the dual standard $t$-structure on $\Ind(\textup{LCoh}(A))$ corresponds to
the $t$-structure defined in Example~\ref{dualstand} (see \cite[X, Definition 3.4.16, Proposition 3.4.18]{DAG} and \cite[Lemma 2.23]{BCN}).
\end{Remark}



\begin{Proposition}
\label{dualcompletion}
Consider the functor $\LMod_{\DD_2(A)}\to \mathsf{D}(\LMod_A)$ induced by the universal pairing $A\otimes_k\DD_2(A)\to k$. This exhibits
$\mathsf{D}(\LMod_A)$ as a left completion of $\LMod_{\DD_2(A)}$.
\end{Proposition}

{\it Proof of Theorem~\ref{intro3} and Proposition~\ref{dualcompletion}.}
We consider the composite
\[
\LMod_{\DD_2(A)}\to \mathsf{D}(\mathsf{D}(\LMod_{\DD_2(A)}))\simeq \mathsf{D}(\Ind(\textup{LCoh}(A)))\to \mathsf{D}(\LMod_A)
\]
where the first arrow is the biduality functor, and the final arrow is induced by $\Upsilon_A:\LMod_A\to \Ind(\textup{LCoh}(A))$.
By Proposition~\ref{pullbackdual},
$\mathsf{D}(\mathsf{D}_{\textup{t}+}(\LMod_{\DD_2(A)}))\simeq \mathsf{D}(\LMod_A)$.
Since $\mathsf{D}(\LMod_A)$ is right complete, $\mathsf{D}_{\textup{t}+}(\mathsf{D}_{\textup{t}+}(\LMod_{\DD_2(A)}))\simeq \mathsf{D}(\LMod_A)$.
We now prove Proposition~\ref{dualcompletion}: We observe that the composite $c:\LMod_{\DD_2(A)}\to \mathsf{D}(\LMod_A)$ is determined by the pairing 
$\LMod_A\otimes\LMod_{\DD_2(A)}\to \Mod_k$ induced by the universal pairing $A\otimes_k\DD_2(A)\to k$. Once this is proven, our assertion follows from 
the fact that $c$ is the left completion of $\LMod_{\DD_2(A)}$ (see \cite[Lemma 3.10]{HD}).
The biduality functor is determined by the universal pairing $V:\LMod_{\DD_2(A)}\otimes\End_{\LMod_{\DD_2(A)}}(\Mod_k)\to \Mod_k$.
By construction, the composite $c$ is induced by the pairing
$V\circ (\textup{id}_{\LMod_{\DD_2(A)}}\otimes (\chi_A\circ \Upsilon_A)):\LMod_{\DD_2(A)}\otimes\LMod_A\to \LMod_k$.
The fully faithfulness of $\chi_A\circ \Upsilon_A$ and $\DD_2(\DD_2(A))\simeq A$ imply that $\chi_A\circ \Upsilon_A$ exhibits $\LMod_A$ as the algebraization of 
$\End_{\LMod_{\DD_2(A)}}(\Mod_k)$ (cf. Definition~\ref{algdef}). It follows from Lemma~\ref{dualcatalg} and $A\simeq \DD_2(\DD_2(A))$ that $c$
is determined by $\LMod_A\otimes\LMod_{\DD_2(A)}\to \Mod_k$ induced by the universal pairing $A\otimes_k\DD_2(A)\to k$.
\QED

\section{Formal stacks}

\label{FS}

In this section, we generalize Theorem~\ref{main0} to categorified sheaves on pointed $\etwo$-formal stacks.
This generalization is intended to include the case of complete algebras; however, our generalization here is rather formal.
We recall the notion of pointed formal $\etwo$-stacks from \cite[X, Section 4]{DAG} where they are referred to as formal moduli problems.
A pointed formal $\etwo$-stack is a functor $F:\textup{Art}_2\to \SSS$ such that 
\begin{itemize}
  \item $F(k)$ is a contractible space,
\item for any pullback square \[
\xymatrix{
R \ar[r] \ar[d] &  R_0 \ar[d] \\
R_1 \ar[r] & R_{01}
}
\]
in $\textup{Art}_2$
such that $H^0(R_0)\to H^0(R_{01})$ and $H^0(R_1)\to H^0(R_{01})$ are surjective, the canonical map
 \[
\xymatrix{
F(R) \ar[r] & F(R_0)\times_{F(R_{01})}F(R_1) 
}
\]
is an equivalence.
\end{itemize}
Let $\widehat{\mathsf{Stack}}_2^\ast$ be the full subcategory of $\Fun(\textup{Art}_2, \SSS)$, which consists of pointed formal $\etwo$-stacks.
For example, for $A\in \Art_2$, the functor $\Spec A$ corepresented by $A$ is a pointed formal $\etwo$-stack so that $(\Art_2)^{op}\subset \widehat{\mathsf{Stack}}_2^\ast$.
By \cite[X, Theorem 4.0.8]{DAG}, there exists an equivalence of $\infty$-categories $\widehat{\mathsf{Stack}}_2^\ast \simeq \Alg_2(\Mod_k)_{/k}$.

\begin{Definition}
We define categorified structure sheaves and categorified $!$-structure sheaves on pointed formal $\etwo$-stacks.
Let $F\in \widehat{\mathsf{Stack}}_2^\ast$ be a pointed formal $\etwo$-stack.
Write $(\Art_2)^{op}_{/F}=(\Art_2)^{op}\times_{\widehat{\mathsf{Stack}}_2^\ast}(\widehat{\mathsf{Stack}}_2^\ast)_{/F}$.
We consider 
\[
((\Art_2)^{op}_{/F})^{op}\stackrel{\textup{forget}}{\longrightarrow} \Art_2\stackrel{\theta^{2}_k}{\longrightarrow} \Alg_1(\PR_k)_{/\Mod_k}.
\]
We define the categorified structure sheaf of $F$ to be the limit of the composite functor. We denote it by $\OO_F^{\textup{cat}}$. 
Similarly, we define $\OO_F^{!,\textup{cat}}$ to be the limit of the composite functor  
\[
((\Art_2)^{op}_{/F})^{op}\stackrel{\textup{forget}}{\longrightarrow} \Art_2\stackrel{\theta^{\Ind,2}_k}{\longrightarrow} \Alg_1(\PR_k)_{/\Mod_k}.
\] 
See Construction~\ref{finalconst} for $\theta^2_k$ and  $\theta^{\Ind,2}_k$.
\end{Definition}

\begin{Remark}
Since $\Alg_1(\PR_k)_{/\Mod_k}$ has small limits, we can take
limits of $\OO_F^{\textup{cat}}$ and $\OO_F^{!,\textup{cat}}$.
The limit $\QC^\otimes(F)$ is the monoidal $\infty$-categories of quasi-coherent sheaves on $F$ (equipped with
the augmentation to $\Mod_k^\otimes$), respectively.  
\end{Remark}

Consider $\mathsf{D}\circ \mathsf{D}:\Alg_1(\PR_k)_{/\Mod_k}\to\Alg_1(\PR_k)_{/\Mod_k}$. It gives rise to
the biduality map $\delta_{\OO_F^{\textup{cat}}}:\OO_F^{\textup{cat}}\to \DDDD(\mathsf{D}(\OO_F^{\textup{cat}}))$.

\begin{Corollary}
\label{mainformal}
There  exists a canonically defined equivalence
$\OO_F^{\textup{cat}}\simeq \DDDD(\mathsf{D}(\OO_F^{\textup{cat}}))$.
\end{Corollary}

\appendix

\section{}


\label{app}




We construct a covariantly functorial assignment $A \mapsto \Ind(\textup{RCoh}(A))$ 
for Artin $\eone$-algebras $A$.
In the commutative case, it is done in \cite{Gai2}.
In this appendix, we indicate how one can adapt methods in \cite{Gai2} to the context of $\eenu$-algebras.
For simplicity and for our purpose, we concentrate on Artin algebras.
In the main text, we use a left version $A \mapsto \Ind(\textup{LCoh}(A))$ obtained
from the right version by passing to opposite algebras.

\begin{Lemma}
\label{dualind}
Let $A,A'$ be Artin $\eone$-algebras.
The following holds.
\label{applem}
\begin{enumerate}
  \item There exists an equivalence 
  \[
  \Ind(\textup{RCoh}(A))\otimes_k\Ind(\textup{RCoh}(A'))\simeq \Ind(\textup{RCoh}(A\otimes_kA'))
  \] given by $p_1^!(-)\otimes p_2^!(-)$
  where $p_1^\sharp:A\otimes_kk\to A\otimes_kA'\leftarrow k\otimes_kA':p_2^\sharp$.

\item Taking $k$-linear duals we have an equivalence $\textup{LCoh}(A)\to \textup{RCoh}(A)^{op}$
defined by $M\mapsto \textup{Hom}_k(M,k)$, where $\textup{Hom}_k(M,k)$ denotes the Hom complex.
This equivalence induces
\[
\Ind(\textup{RCoh}(A))^\vee=\Ind(\textup{RCoh}(A)^{op})\simeq \Ind(\textup{LCoh}(A)).
\]
Namely, the dual of $\Ind(\textup{RCoh}(A))$ in $\PR_k$ is equivalent to $\Ind(\textup{LCoh}(A))$.

\item For a morphism $f^{\sharp}:A\to B$ of Artin $\eone$-algebras,
the dual 
\[
f_*^\vee:\Ind(\textup{RCoh}(A))\simeq \Ind(\textup{LCoh}(A))^\vee\to \Ind(\textup{LCoh}(B))^\vee\simeq \Ind(\textup{RCoh}(B))
\]
of $f_*:\Ind(\textup{LCoh}(B))\to \Ind(\textup{LCoh}(A))$
is the (right module version of) $!$-pullback functor, that is, the right adjoint of
the (right module version of) the pushforward functor $\Ind(\textup{RCoh}(B))\to \Ind(\textup{RCoh}(A))$.
  
\end{enumerate}

\end{Lemma}

\Proof
We can prove the assertion (1) in the same way as \cite[Vol.1, Chapter 4, Proposition 6.3.4]{Gai2}.
We prove (2). Note first that compactly generated stable $\infty$-categories are dualizable,
and $\Ind(\textup{RCoh}(A))^\vee=\Fun_k(\Ind(\textup{RCoh}(A)),\Mod_k)$
is equivalent to $\Ind(\textup{RCoh}(A)^{op})$.
It is enough to show that $\textup{RCoh}(A)^{op}\simeq \textup{LCoh}(A)$.
We have an equivalence $\textup{LCoh}(A)\to \textup{RCoh}(A)^{op}$
defined by $M\mapsto \textup{Hom}_k(M,k)$.
Finally, we remark that $\textup{Hom}_k(M,k)$ is regarded as a sort of Grothendieck duality:
Since the $!$-pullback $\Mod_k\to\Ind(\textup{RCoh}(A))$ along $A\to k$ is a right adjoint of the Ind-extension $\Ind(\textup{RCoh}(A))\to \Mod_k$
of the forgetful functor $\textup{RCoh}(A)\to \Mod_k$,
thus the image of $k$ under the $!$-pullback is $\omega_A\simeq \Hom_{k}(A,k)$ (note also that $\Hom_{k}(A,k)$ has a diagonal $A$-$A$-bimodule structure). Then we have equivalences 
\begin{eqnarray*}
\operatorname{Hom}_{\LMod_A}(M, \omega_A) &\simeq& \operatorname{Hom}_{\LMod_A}(M, \operatorname{Hom}_k(A, k)) \\ 
&\simeq& \operatorname{Hom}_k(A \otimes_A M, k) \\ &\simeq& \operatorname{Hom}_k(M, k).
\end{eqnarray*}
See \cite[9.3.3]{Ind} for (3); it holds in our context, using the duality in (2).
\QED

\begin{Construction}
\label{enconst}
For $A\in \Alg_1(\Mod_k)$,
we define $\Psi_A:\Ind(\textup{LCoh}(A))\to \LMod_A$ to be the (colimit-preserving) $\Mod_k$-linear functor that extends the inclusion $\textup{LCoh}(A)\hookrightarrow \LMod_A$.
Note that $\Ind(\textup{LCoh}(A))$ is dualizable in $\PR_k$, with dual $\Ind(\textup{LCoh}(A))^\vee$ given by 
$\Ind(\textup{RCoh}(A))$ (see Lemma~\ref{applem} (2)).
We also note that $\LMod_A$ is dualizable in $\PR_k$, with dual $\LMod_A^\vee$ given by 
$\RMod_A$.
Let $\textup{Pr}_k^{\textup{L}, \textup{dual}}$ denote the full subcategory of $\PR_k$, which consists of
dualizable objects.
Let $(-)^\vee:(\textup{Pr}_k^{\textup{L}, \textup{dual}})^{op}\simeq \textup{Pr}_k^{\textup{L}, \textup{dual}}$
be the equivalence given by taking duals.
Taking duals $(-)^\vee$,
we define $\Upsilon_A:\RMod_A\to \Ind(\textup{RCoh}(A))$ to be
\[
\Psi_A^\vee:\RMod_A\simeq \LMod_A^\vee\to  \Ind(\textup{LCoh}(A))^\vee\simeq \Ind(\textup{RCoh}(A)).
\]
When $n=\infty$ and $A\in \Art=\Art_n$,
\[
\Upsilon_A:\Mod_A\longrightarrow \Ind(\Coh(A))
\]
coincides with $\Upsilon_A$ in \cite{Gai2} given by $M\mapsto M\otimes_A\omega_A$.
\end{Construction}

\begin{Proposition}
\label{pullbackdual}
Suppose that $A$ is an Artin $\eone$-algebra.
Then the functor $\Upsilon_A:\RMod_A\to \Ind(\textup{RCoh}(A))$
is fully faithful which carries $A$ to $\omega_A$. In particular, $\Upsilon_{A}$ gives an equivalence
$\RMod_A^{\le0}\simeq \Ind(\textup{RCoh}(A))^{\le0,ds}$ (see Definition~\ref{dualtind})
so that $\Upsilon_A$ exhibits $\RMod_A$ as a right completion of $\Ind(\textup{RCoh}(A))^{ds}$.
\end{Proposition}

\Proof
In fact, $A\in \RMod_A$ corresponds to the forgetful functor $\LMod_A\to \Mod_k$
through $\RMod_A\simeq \Fun_k(\LMod_A,\Mod_k)$.
The composition with $\Psi_A$ induces
\[
\Fun_k(\LMod_A,\Mod_k)\to \Fun_k(\Ind(\textup{LCoh}(A)),\Mod_k)
\]
which carries the forgetful functor $\LMod_A\to \Mod_k$ to 
the Ind-extension $h:\Ind(\textup{LCoh}(A)) \to \Mod_k$ which extends
the forgetful functor $\textup{LCoh}(A)\to \Mod_k$.
Thorugh Yoneda embedding
$\textup{LCoh}(A)^{op}\subset \Fun_k(\Ind(\textup{LCoh}(A)),\Mod_k)$,
$h$ corresponds to $A\in \textup{LCoh}(A)^{op}$.
Through $\textup{LCoh}(A)^{op}\simeq \textup{RCoh}(A)\subset \Ind(\textup{RCoh}(A))$,
$A$ maps to $\Hom_k(A,k)\simeq \omega_A$ (see Lemma~\ref{dualind} (2) and its proof).
We observe $\Upsilon_A$ is fully faithful.
To see this, we first note that the colimit-preserving functor $\Xi_A:\LMod_A\to \Ind(\textup{LCoh}(A))$
extending the inclusion $\LMod_A^{\omega}=\textup{LPerf}_A\hookrightarrow \textup{LCoh}(A)$ is a left adjoint to $\Psi_A$.
By definition, $\Psi_A\circ \Xi_A$ is equivalent to the identity functor.
Taking the duals $\Xi_A^\vee:\Ind(\textup{RCoh}(A))\to \RMod_A$ is a right adjoint to $\Upsilon_A$.
Thus, $(\Upsilon_A,\Xi_A^\vee)$ is a colocalization adjoint pair.
In particular, $\Upsilon_A$ is fully faithful.
\QED


We will give a construction of a covariantly functorial assignment $A\mapsto \Ind(\textup{RCoh}(A))$ (for $A\in \textup{Art}_1$) as a symmetric monoidal functor
such that
the covariant functoriality is given by $!$-pullback functors
(see Construction~\ref{finalconst}).
This requires some preliminaries.

\begin{Lemma}
\label{keyapp}
Let $A_1,\ldots, A_m,B$ be Artin $\eone$-algebras.
Let 
\[
f:\otimes_{1\le i\le m}\RMod_{A_i}=\RMod_{A_1}\otimes_k \cdots  \otimes_k\RMod_{A_m} \to \RMod_B
\]
be a morphism in $(\PRT)_k$.
Suppose that the dual (in $\PR_k$) $f^\vee:\LMod_B\to \LMod_{A_1\otimes_k\cdots \otimes_kA_m}\simeq\LMod_{A_1}\otimes_k \cdots  \otimes_k\LMod_{A_m}$
is left $t$-exact up to a finite shift with respect to the standard $t$-structures.
Then there exists a morphism $\overline{f}:\Ind (\textup{RCoh}(A_1))\otimes_k \cdots  \otimes_k\Ind(\textup{RCoh}(A_m)) \to  \Ind(\textup{RCoh}(B))$ in $\PR_k$
such that
\[
\xymatrix{
\RMod_{A_1}\otimes_k \cdots  \otimes_k \RMod_{A_m}   \ar[d]^{\Upsilon_{A_1}\otimes \cdots \otimes\Upsilon_{A_m}}  \ar[r]^{f} & \RMod_B \ar[d]^{\Upsilon_B} \\
\Ind (\textup{RCoh}(A_1))\otimes_k \cdots  \otimes_k\Ind(\textup{RCoh}(A_m)) \ar[r]^(0.7){\overline{f}} &  \Ind(\textup{RCoh}(B))
}
\]
commutes.
Such $\overline{f}$ is unique up to a contractible space of choices.
\end{Lemma}

\Proof
Taking duals,
we consider the diagram
\[
\xymatrix{
\LMod_{A_1}\otimes_k \cdots  \otimes_k \LMod_{A_m}    &  \LMod_B \ar[l]^{f^\vee} \\
\Ind (\textup{LCoh}(A_1))\otimes_k \cdots  \otimes_k\Ind(\textup{LCoh}(A_m)) \ar[u]^{\Psi_{A_1}\otimes \cdots \otimes\Psi_{A_m}}   &  \Ind(\textup{LCoh}(B)) \ar[u]^{\Psi_B} 
}
\]
The left vertical arrow can be identified with
\[
\Psi_{A_1\otimes_k \cdots  \otimes_kA_m}:\Ind (\textup{LCoh}(A_1\otimes_k \cdots  \otimes_kA_m))  \to \LMod_{A_1\otimes_k \cdots  \otimes_kA_m}
\]
through canonical equivalences $\Ind (\textup{LCoh}(A_1\otimes_k \cdots  \otimes_kA_m))\simeq \Ind (\textup{LCoh}(A_1))\otimes_k \cdots  \otimes_k\Ind(\textup{LCoh}(A_m))$
(see Lemma~\ref{dualind}) and $\LMod_{A_1\otimes_k\cdots \otimes_kA_m}\simeq\LMod_{A_1}\otimes_k \cdots  \otimes_k\LMod_{A_m}$.
The functor $\Psi_{A_1\otimes_k \cdots  \otimes_kA_m}$ induces an equivalence of full subcategories of left bounded objects
$\Ind (\textup{LCoh}(A_1\otimes_k \cdots  \otimes_kA_m))^+ \to \LMod_{A_1\otimes_k \cdots  \otimes_kA_m}^+$
with respect to the standard $t$-structures.
By our assumption, $f^\vee$ sends $\LMod_B^+$ to  $\LMod_{A_1\otimes_k \cdots  \otimes_kA_m}^+$.
In addition, $\textup{LCoh}_B$ maps to $\LMod_B^+$.
Then we have an essentially unique $f'$ which makes the diagram 
\[
\xymatrix{
   \LMod_{A_1\otimes_k \cdots  \otimes_kA_m}^+  &  \LMod_B^+ \ar[l]^{f^\vee} \\
\Ind (\textup{LCoh}(A_1\otimes_k \cdots  \otimes_kA_m))^+ \ar[u]^{\Psi_{A_1\otimes_k \cdots  \otimes_kA_m}}_{\simeq}   &  \textup{LCoh}(B) \ar[u]^{\Psi_B} \ar[l]^(0.3){f'} 
}
\]
commute.
Then we take $\overline{f}$ to be the dual of $\Ind(\textup{LCoh}(B))\to \Ind (\textup{LCoh}(A_1\otimes_k \cdots  \otimes_kA_m))$ induced by $f'$, which has the required uniqueness by construction.
\QED

Let $(\PRTT)_{k,\textup{Art}}$ be the full subcategory of $(\PRTT)_k$ which consists of objects equivalent to the form $(\RMod_A,\RMod_A^{\le0})$ with $A\in \textup{Art}_1$.
Here $\RMod_A^{\le0}$ is the connective part of the standard $t$-structure (see Example~\ref{connectivetstr}).
This full subcategory $(\PRTT)_{k,\textup{Art}}$ is closed under taking tensor product and contains the unit object $\Mod_k$ with the standard $t$-structure.
Thus, it inherits a symmetric monoidal structure.
We write $(\PRTT)_{k,\textup{Art}}^\otimes\to \Gamma$ for the corresponding coCartesian fibration (see Section~\ref{notation} for $\Gamma$).

Let $(\PRT)_{k,\textup{Art},\Ind}$ be the (non-full) subcategory of $(\PRT)_k$, which consists of objects equivalent to $\Ind(\textup{RCoh}(A))^{ds}=(\Ind(\textup{RCoh}(A)),\Ind(\textup{RCoh}(A))^{\le0,ds})$ with some $A\in \textup{Art}_1$.
Here $\Ind(\textup{RCoh}_A)^{\le0,ds}$ is the connective part of
the {\it dual} standard $t$-structure (see Definition~\ref{dualtind}).
The space of morphisms 
\[
\Map_{(\PRT)_{k,\textup{Art},\Ind}}(\Ind(\textup{RCoh}(A))^{ds},\Ind(\textup{RCoh}(B))^{ds})
\]
is
the full subcategory of $\Map_{(\PRT)_{k}}(\Ind(\textup{RCoh}(A))^{ds},\Ind(\textup{RCoh}(B))^{ds})$
consisting of right $t$-exact functors
corresponding to those ($1$-)morphisms $f$ satisfying the following conditions:
\begin{itemize}
\item
 $f^\vee:\Ind(\textup{LCoh}(B))\simeq \Ind(\textup{RCoh}(B))^\vee \to  \Ind(\textup{RCoh}(A))^\vee \simeq  \Ind(\textup{LCoh}(A))$  preserves compact objects,

\item $f^\vee$ is left $t$-exact up to a finite shift with respect to the standard $t$-structure.

\end{itemize}

\begin{Claim}
The subcategory $(\PRT)_{k,\textup{Art},\Ind}$ inherits a symmetric monoidal structure from $(\PRT)_{k}$.
\end{Claim}

\Proof
Since $\Ind(\textup{RCoh}(A))\otimes_k \Ind(\textup{RCoh}(B))\simeq \Ind(\textup{RCoh}(A\otimes_kB))$ (see \cite{Ind}, Lemma~\ref{dualind} (1)),
objects of $(\PRT)_{k,\textup{Art},\Ind}$ are closed under tensor products in $(\PRT)_{k}$, and the unit $(\Mod_k,\Mod_k^{\le0})$ belongs to $(\PRT)_{k,\textup{Art},\Ind}$.
It remains to verify that for any object $\Ind(\textup{RCoh}(C))\in (\PRT)_{k,\textup{Art},\Ind}$
and any morphism $f:\Ind(\textup{RCoh}(A))\to \Ind(\textup{RCoh}(B))$ in $(\PRT)_{k,\textup{Art},\Ind}$,
$\textup{id}_{\Ind(\textup{RCoh}(C))}\otimes f$ is also a morphism in $(\PRT)_{k,\textup{Art},\Ind}$.
We consider
\[
\textup{id}_{\Ind(\textup{LCoh}(C))}\otimes f^\vee:\Ind(\textup{LCoh}(C))\otimes_k\Ind(\textup{LCoh}(B))\to \Ind(\textup{LCoh}(C))\otimes_k\Ind(\textup{LCoh}(A)).
\]
The restricted functor 
\[
\textup{id}\otimes f^\vee:\textup{LCoh}(C)\otimes_k\textup{LCoh}(B)\to \textup{LCoh}(C)\otimes_k\textup{LCoh}(A)
\]
carries objects lying in the amplitude $[r,\infty)$ to objects lying in the amplitude $[r-N,\infty)$ for some $N>>0$. (Here the tensor product means the tensor product of small stable idempotent-complete $\infty$-categories over $k$.)
The canonical functors
$\textup{LCoh}(C)\otimes_k\textup{LCoh}(A)\to \textup{LCoh}(C\otimes_kA)$ and 
$\textup{LCoh}(C)\otimes_k\textup{LCoh}(B)\to \textup{LCoh}(C\otimes_kB)$ are equivalences
(cf. e.g. the proof of \cite[Vol.1, Chapter 4, Proposition 6.3.4]{Gai2}).
Also, $\Ind(\textup{LCoh}(D))^{\ge r}=\Ind(\textup{LCoh}(D)^{\ge r})$ for any $D\in \textup{Art}_1$
when it is equipped with the standard $t$-structure.
Thus, we see that $\textup{id}\otimes f^\vee$ is left $t$-exact up to a shift $[N]$.
\QED

We write
$(\PRT)_{k,\textup{Art},\Ind}^\otimes\to \Gamma$
for the coCartesian fibration of the resulting symmetric monoidal $\infty$-category.

Let 
\[
R:(\PRT)_{k,\textup{Art},\Ind}^\otimes\longrightarrow  (\PRTT)_{k,\textup{Art}}^\otimes
\]
be the map over $\Gamma$, which is the lax symmetric monoidal functor
induced by the lax symmetric monoidal right completion functor $\widehat{R}:(\PRT)^\otimes\to (\PRTT)^\otimes$ (cf. Section~\ref{tsec}).
This lax symmetric monoidal functor $R$ is symmetric monoidal.
Indeed, the composite of canonical functors
\begin{eqnarray*}
\LMod_A\otimes_k\LMod_B &\simeq& \widehat{R}(\Ind(\textup{LCoh}(A)))\otimes_k\widehat{R}(\Ind(\textup{LCoh}(B))) \\
 &\to& \widehat{R}(\Ind(\textup{LCoh}(A))\otimes_k\Ind(\textup{LCoh}(B))) \\
&\simeq& \widehat{R}(\Ind(\textup{LCoh}(A\otimes_kB)))  \\
&\simeq& \LMod_{A\otimes_kB}
\end{eqnarray*}
is an equivalence
(observe that $A\otimes_kB\in \LMod_{A\otimes_kB}=\LMod_A\otimes_k\LMod_B$ maps to $A\otimes_kB$ which induces the identity of the endomorphism algebras).

\begin{Lemma}
\label{fullyfaithfullemma}
The map
$R:(\PRT)_{k,\textup{Art},\Ind}^\otimes\longrightarrow  (\PRTT)_{k,\textup{Art}}^\otimes$ is $1$-fully faithful
as a functor of $\infty$-categories. Namely, for any $C,C'\in(\PRT)_{k,\textup{Art},\Ind}^\otimes$,
the induced map $\Map_{(\PRT)_{k,\textup{Art},\Ind}^\otimes}(C,C')\to  \Map_{(\PRTT)_{k,\textup{Art}}^\otimes}(R(C),R(C'))$
is a fully faithful functor between mapping spaces.
\end{Lemma}

\Proof
We first note that by definition $\Upsilon_A$ is 
$\Psi_A^\vee:\RMod_A=\LMod_A^\vee\to \Ind(\textup{LCoh}(A))^{\vee}=\Ind(\textup{RCoh}(A))$,
and $\Upsilon_A:\RMod_A\to \Ind(\textup{RCoh}(A))^{ds}$ exhibits $\RMod_A$ as a right completion, where the superscript $ds$ indicates the dual standard $t$-structure (see Definition~\ref{dualtind}, Proposition~\ref{pullbackdual}).
Let $A_1,\ldots, A_m,B$ be Artin $\eone$-algebras.
If we consider the tensor product
 $\otimes_{1\le i\le m}\Upsilon_{A_i}:\otimes_{1\le i\le m}\RMod_{A_i} \to \otimes_{1\le i\le m}\Ind(\textup{RCoh}(A_i))$,
then, either by the dual of [Gai2, Vol.1, Chapter 4, Proposition 6.3.4] or by regarding each $\Upsilon_{A_i}$ as a $k$-linear functor carrying $A_i \in \RMod_{A_i}$ to $\omega_{A_i}$ and inducing an equivalence on endomorphism algebras, this tensor product can naturally be identified with $\Upsilon_{A_1\otimes_k\cdots\otimes_kA_m}:\RMod_{A_1\otimes_k\cdots\otimes_kA_m}\to \Ind(\textup{RCoh}(A_1\otimes_k\cdots\otimes_kA_m))$
 through the equivalences $\otimes_{1\le i\le m}\RMod_{A_i}\simeq \RMod_{A_1\otimes_k\cdots\otimes_kA_m}$
and $\otimes_{1\le i\le m}\Ind(\textup{RCoh}(A_i)) \simeq \Ind(\textup{RCoh}(A_1\otimes_k\cdots\otimes_kA_m))$.
The functor $R$ induces a map of mapping spaces
\begin{eqnarray*}
R_{A_1,\ldots,A_m;B}:\Map_{(\PRT)_{k,\textup{Art},\Ind}}(\Ind (\textup{RCoh}(A_1))^{ds}\otimes_k \cdots  \otimes_k\Ind(\textup{RCoh}(A_m))^{ds},\Ind(\textup{RCoh}(B))^{ds}) \to \ \ \ \ \ \ \ \ \  \ \   \\
\ \ \ \ \ \ \   \Map_{(\PRTT)_{k,\textup{Art}}}(\RMod_{A_1}\otimes_k \cdots  \otimes_k \RMod_{A_m}, \RMod_B).
\end{eqnarray*}
To prove our assertion, it is enough to show that $R_{A_1,\ldots,A_m;B}$ is fully faithful when we regard the mapping spaces as $\infty$-categories.
The essential image is contained in the full subcategory of the target spanned by those functors
$\RMod_{A_1}\otimes_k \cdots  \otimes_k \RMod_{A_m}\simeq \RMod_{A_1\otimes_k\cdots \otimes_kA_m}\to \RMod_B$ which are left $t$-exact up to finite shifts.
Indeed, 
for $A\in \textup{Art}_1$,
$\Psi_A:\Ind(\textup{LCoh}(A))^{st}\to \LMod_A$ is a left completion, where the superscript $st$ indicates the standard $t$-structure.
Thus, if a $k$-linear functor $f:\Ind(\textup{LCoh}(A))^{st}\to \Ind(\textup{LCoh}(B))^{st}$ is left $t$-exact up to a finite shift,
the induced functor $\widehat{R}(f):\LMod_A\to \LMod_B$ is left $t$-exact up to a finite shift.
From Lemma~\ref{keyapp} we conclude that $R_{A_1,\ldots,A_m;B}$ is fully faithful.
\QED

\begin{Construction}
\label{thetaTTT}
We will construct
$\Theta_k^{+}:\Alg_{1}(\Mod_k)\to  \Mod_{\Mod_k}(\PRTT)_{\Mod_k/}$,
which carries $B$ to $(\RMod_B, \RMod_B^{\le0})$ together with $(-)\otimes_kB:\Mod_k\to \RMod_B$,
where $\RMod_B^{\le0}$ is the smallest full subcategory which contains $B$ and is closed under small colimits and extensions (see Remark~\ref{unifiedtstr}).
Recall that
$\Theta_k^{\textup{right}}=\Theta_k\circ \operatorname{opp}:\Alg_{1}(\Mod_k)\stackrel{\operatorname{opp}}{\simeq} \Alg_1(\Mod_k)\to  \Mod_{\Mod_k}(\PR)_{\Mod_k/}= (\PR_k)_{\Mod_k/}$
is fully faithful where the equivalence $\Alg_{1}(\Mod_k)\simeq \Alg_1(\Mod_k)$
is determined by taking opposite algebras. We denote the essential image by $T$. Let $S$ be the full subcategory of $\Mod_{\Mod_k}(\PRTT)_{\Mod_k/}$
spanned by objects of the form $(\RMod_B, \RMod_B^{\le0})$ with $B\in\Alg_{1}(\Mod_k)$.
The forgetful functor $\Mod_{\Mod_k}(\PRTT)_{\Mod_k/}\to (\PR_k)_{\Mod_k/}$ induces $f:S\to T$.
We claim that $f$ is an equivalence. 
We define $\Theta_k^{+}$ to be the composition $g\circ \Theta_k^{\textup{right}}|_{\Alg_{1}(\Mod_k)}$,
where $g$ is an inverse functor of $f$.
Clearly, $f$ is essentially surjective.
Note that every morphism $F:\RMod_A\to \RMod_B$ in $(\PR_k)_{\Mod_k/}$
sends $A$ to $B$. Since $\RMod_A^{\le0}$ (resp. $\RMod_B^{\le0}$)
is the smallest full subcategory of $\RMod_A$ (resp. $\RMod_B$)
which contains $A$ (resp. $B$) and is closed under small colimits and extensions,
it follows that $F$ is right $t$-exact. Thus, $f$ is fully faithful.

Since the forgetful functor $\Mod_{\Mod_k}(\PRTT)_{\Mod_k/}\to \Mod_{\Mod_k}(\PR)_{\Mod_k/}$ is (promoted to) a symmetric monoidal functor,
and $\Theta_k^{\textup{right}}$ is promoted to a symmetric monoidal functor (see \cite[4.8.5.16]{HA}),
thus $\Theta_k^{+}$ is promoted to a symmetric monoidal functor $\Alg_1^+(\Mod_k)^\otimes\to (\Mod_{\Mod_k}(\PRTT)_{\Mod_k/})^\otimes$.
\end{Construction}

\begin{Construction}
\label{finalconst}
Let $\textup{Art}_1^\otimes$ denote the symmetric monoidal full subcategory of $\Alg_1^+(\Mod_k)^\otimes$, spanned by Artin $\eone$-algebras. (Note the equivalence $\Art_n\simeq \Art_n^\circ$ (cf. Remark~\ref{augnonaug}).)
We regard an object of $\textup{Art}_1^\otimes$ lying over $\langle m\rangle\in \Gamma$ as $(A_1,\ldots,A_m,\langle m\rangle)$.
By Construction~\ref{thetaTTT}, we have a symmetric monoidal functor 
\[
\theta_k^+=\Theta_k^{+}|_{\textup{Art}_1^\otimes}:\textup{Art}_1^\otimes\to (\PRTT)_{k,\textup{Art}}^\otimes,
\]
which is defined as a map of coCartesian fibrations over $\Gamma$.

We consider the assignment $F:\textup{Art}_1^+(\Mod_k)^\otimes\to (\PRT)_{k,\textup{Art},\Ind}^\otimes$  at the level of objects given by
\[
\textup{Art}_1^+(\Mod_k)^\otimes\ni (A_1,\ldots,A_m,\langle m\rangle)\mapsto (\Ind(\textup{RCoh}(A_1))^{ds},\ldots, \Ind(\textup{RCoh}(A_m))^{ds},\langle m\rangle)\in(\PRT)_{k,\textup{Art},\Ind}^\otimes.
\]
Next, we define the assignment $F$ at the level of $1$-morphisms.
To a morphism $(A_1,\ldots,A_m,\langle m\rangle)\to (B_1,\ldots,B_n,\langle n\rangle)$ over $\alpha:\langle m\rangle\to \langle n\rangle$
we associate the morphism determined by the $!$-pullback
\[
P_i:\otimes_{j\in \alpha^{-1}(i)}\Ind(\textup{RCoh}(A_{j}))^{ds}\to \Ind(\textup{RCoh}(B_i))^{ds}
\]
along $\otimes_{j\in \alpha^{-1}(i)}A_{j}\to B_i$ for each $i\in  \langle n\rangle^\circ$. 
Notice that the dual of $P_i$ is the pushforward functor (i.e., the Ind-extension of the restriction along $\otimes_{j\in \alpha^{-1}(i)}A_{j}\to B_i$)
so that the associated morphism lies in $(\PRT)_{k,\textup{Art},\Ind}^\otimes$.
The assignment obtained by applying the right completion $\widehat{R}:(\PRT)_{k}^\otimes\to (\PRTT)_{k}^\otimes$ to $F$ is extended to $\theta_k^{+}$. By the fact from Lemma~\ref{fullyfaithfullemma} that $R$ is $1$-fully faithful (cf. \cite[Vol.1 Chapter 4, Lemma 2.2.7]{Gai2}),
the assignment $F$ is extended to a symmetric monoidal functor $\textup{Art}_1^\otimes\to(\PRT)_{k,\textup{Art},\Ind}^\otimes$.
Composing with the forgetful functor $(\PRT)_{k,\textup{Art},\Ind}^\otimes\to (\PRT)_k^\otimes$, we obtain
\[
\theta_k^{\Ind}:\textup{Art}_1^\otimes \longrightarrow (\PRT)_k^\otimes.
\]
By construction, it carries $A$ to $\Ind(\textup{RCoh}(A))^{ds}$ and sends a morphism $A\to B$ to the $!$-pullback functor $\Ind(\textup{RCoh}(A))\to \Ind(\textup{RCoh}(B))$.
Using the right completion,
$(\PRT)_k^\otimes\stackrel{\widehat{R}}{\longrightarrow} (\PRTT)_k^\otimes\hookrightarrow (\PRT)_k^\otimes$, 
we have a natural transformation $\Upsilon$ from $\theta_k:\textup{Art}_1^\otimes \longrightarrow (\PRT)_k^\otimes$ (we abuse notation)
to $\theta_k^{\Ind}$. For any $A\in \textup{Art}_1$, the evaluation of $\Upsilon$ at $A$ is equivalent to $\Upsilon_A$ defined in Construction~\ref{enconst}. 
\end{Construction}

\begin{Remark}
Using the symmetric monoidal equivalence $\textup{Art}_1^\otimes\simeq \textup{Art}_1^\otimes$ given by passing to opposite algebras $A\mapsto A^{op}$,
we obtain left module versions of $\theta_k^{\Ind}$, $\theta_k$ and $\Upsilon$. 
In the main text, we abuse notation by writing $\theta_k^{\Ind}$, $\theta_k$ and $\Upsilon$
for  left module versions of $\theta_k^{\Ind}$, $\theta_k$ and $\Upsilon$, respectively.
\end{Remark}

\begin{Definition}
\label{eenuindmodule}
Let
\[
\theta_k^{n}:\Art_n=\Alg_{n-1}(\textup{Art}_1^\otimes) \to \Alg_{n-1}((\PR)_k^\otimes)
\]
be the functor induced by $\theta_k^+$ and the forgetful functor (see Construction~\ref{finalconst}),
which carries $A$ to $\RMod_A$.
We define
\[
\theta_k^{\Ind,n}:\Art_n=\Alg_{n-1}(\textup{Art}_1^\otimes) \to \Alg_{n-1}((\PR)_k^\otimes)
\]
to be the functor
induced by $\theta_k^{\Ind}$ and the forgetful functor (see Construction~\ref{finalconst}),
which carries $A$ to $\Ind(\textup{RCoh}(A))$.
We abuse notation and write
\[
\Upsilon:\theta_k^{n}\longrightarrow \theta_k^{\Ind,n}
\]
 the natural transformation
induced by $\Upsilon$ (see Construction~\ref{finalconst}).
\end{Definition}

\begin{Definition}
For $A\in \Art_n$, we define the $\mathbf{E}_{n-1}$-monoidal $\infty$-category $\Ind(\textup{RCoh}(A))$ to be $\theta_k^{\Ind,n}(A) \in \Alg_{n-1}((\PR)_k^\otimes)$ in Definition~\ref{eenuindmodule}.
\end{Definition}

\begin{Remark}
Recall the equivalence $\Art_n\simeq \Art_n^\circ$ (cf. Remark~\ref{augnonaug}).
Using this equivalence, one can promote $\theta_k^n$ to $\Art_n\to \Alg_{n-1}(\textup{Pr}_k^{\textup{L}, \textup{dual}})_{/\Mod_k}$.
Similarly, 
$\theta_k^{\Ind,n}$ is naturally extended to a functor 
\[
\Art_n\to \Alg_{n-1}(\textup{Pr}_k^{\textup{L}, \textup{dual}})_{/\Mod_k}.
\]
\end{Remark}

\


\vspace{2mm}

{\it Acknowledgements.}
The author would like to thank Takuo Matsuoka. This work originated from email correspondence with him.
The author is deeply grateful to the anonymous referee for pointing out a serious error in the original proof of Theorem 1.3 and for providing the elegant argument presented in this paper. The referee's insight and generosity were instrumental in bringing this paper to its present form.
The responsibility for the contents of this paper is entirely the author's.
This work was partially supported by JSPS KAKENHI Grant.

\end{document}